\documentclass[11pt]{amsart}

\usepackage{tikz}
\usepackage{pgfmath-xfp}
\usepackage{pgfplots}
\usetikzlibrary{shadows}
\usetikzlibrary{shapes}
\usetikzlibrary{decorations}
\tikzset{               
    redarrows/.style={postaction={decorate},decoration={markings,mark=at position 0.1 with {\arrow[draw=black]{>}}},
           decoration={markings,mark=at position -0.4 with {\arrow[draw=black]{>}}},}}

\usepackage [utf8]{inputenc}
\usepackage {tikz}
\usepackage{tkz-euclide}
\usepackage[utf8]{inputenc}

\usepackage{xcolor}
\usepackage{cite}
 \usepackage[colorlinks=true]{hyperref}
 \hypersetup{urlcolor=blue, citecolor=red}

\usepackage{amssymb}
\usepackage{amsmath}
\usepackage{amscd}
\usepackage{curves}
\usepackage{epsfig}
\usepackage{bbm}
\usepackage{geometry}
\usepackage{graphicx}
\usepackage{pstricks}
\usepackage{pstricks-add}
\usepackage{pst-grad}
\usepackage{pst-plot}
\usepackage{verbatim}
\usepackage{latexsym}
\usepackage{eucal}
\usepackage{amsfonts}
\usepackage{enumerate}

\numberwithin{equation}{section}
\newtheorem{theorem}{Theorem}[section]
\newtheorem{lemma}[theorem]{Lemma}
\newtheorem{prop}[theorem]{Proposition}

\def \bpf {\begin{proof}}
\def \epf {\end{proof}}
\def \beq {\begin{equation}}
\def \eeq {\end{equation}}
\def \bsp{\begin{split}}
\def \esp{\end{split}}

\def \defn {:=}

\def \wt {\widetilde}

\def \tr {\operatorname{Tr}}

\def \id {{\operatorname{I}}}

\def \mca {{\mathcal A}}
\def \mcb {{\mathcal B}}

\def \mcd {{\mathcal D}}
\def \mce {{\mathcal E}}

\def \mcf {{\mathcal F}}

\def \mch {{\mathcal H}}

\def \mcl {{\mathcal L}}

\def \mcw {{\mathcal W}}

\def \mcp {{\mathcal P}}

\def \mcr {{\mathcal R}}
\def \mcs {{\mathcal S}}
\def \mct {{\mathcal T}}
\def \mcu {{\mathcal U}}

\def \mcz {{\mathcal Z}}

\def \CC {{\mathbb C}}
\def \mh {{\mathbb H}}
\def \mb {{\mathbb B}}
\def \mbb {{\mathbb B}^{n+1}}
\def \mbh {{\mathbb H}}

\def \mr {{\mathbb R}}
\def \mn {{\mathbb N}}
\def \ms {{\mathbb S}}

\def\ha {\frac{1}{2}}

\def \nsq {\frac{n^2}{4}}

\def \ka {\kappa}

\def \Tau {{\mathcal{T}}}

\def \spec {\operatorname{spec}}

\def \im {\operatorname{Im}}
\def \re {\operatorname{Re}}

\def \supp {\text{supp }}

\def \mrn {{\mathbb R}^n}
\def \mhn {{\mathbb H}^{n+1}}

\def \msn {{\mathbb S}^n}

\def \ga {{\gamma}}
\def \eps {\varepsilon}

\def \vphi {\varphi}
\def \vrho {\varrho}

\def \la {\lambda}   
\def \La {\Lambda}   
\def \lan {\langle}   
\def \ran {\rangle}   
\def \del {\delta}   

\def \p {\partial}

\def \novt {\frac{n}{2}}

\def \beqq {\begin{equation}}
\def \bush {\mcb_0}
\def \eeqq {\end{equation}}

\def \mck {\mathcal{K}}

\def \ooh {\frac{1}{h}}

\def \ioh {\frac{i}{h}}

\numberwithin{equation}{section}

\begin{document}

\title[ Phase Shifts in Potential Scattering in $\mbh^n$] {The High Energy Distribution of Scattering Phase Shifts of Schr\"odinger operators in Hyperbolic Space}
\author[S\'a Barreto]{Ant\^onio S\'a Barreto}
\address{Department of Mathematics, Purdue University \newline
\indent 150 North University Street, West Lafayette IN  47907, USA}
\email{sabarre@math.purdue.edu}

\begin{abstract} We prove  a trace formula for the high energy limit of the  scattering phase shifts of Schr\"odinger operators with short range real valued potentials in hyperbolic space; it relates the scattering shifts and the geodesic $X$-ray transform of the potential. This extends a result of  Bulger and Pushnitski \cite{BulPus1}  for Schr\"odinger operators in Euclidean space. As an application, we prove that the high energy limit of the phase shifts uniquely determines radial potentials which are monotone and decay super-exponentially.   This extends a result of Levinson \cite{Lev} for potential perturbations of the Euclidean Laplacian  to this special class of potentials in hyperbolic space.
\end{abstract}
\maketitle

\tableofcontents

\section{Introduction}

When a wave interacts with a perturbation it undergoes a phase shift and this phenomenon has been well studied by physicists and mathematicians, see for example \cite{Bar,Ber,Eu,Landau,Lev,LipS}, largely for perturbations of the Euclidean Laplacian by real valued radially symmetric potentials, but also for more general perturbations \cite{BirYaf2,BulPus1,GelHas0,GelHasZel,GelIng,Ing} of the Euclidean space, and in other settings as well \cite{ChrUri,KosNil}. The inverse problem of determining a potential from its phase shifts has also been considered in \cite{Bar,Lev}.

The phase shifts  are defined to be the logarithm of the eigenvalues of the (relative) scattering matrix,  which is a unitary operator. So  this is also a spectral problem concerning the distribution of the eigenvalues of a physically meaningful operator which is not self adjoint.  Our main result, Theorem \ref{main}, establishes  a classical-quantum trace formula for Schr\"odinger operators with short range potentials in hyperbolic space.  It relates the high energy distribution of the phase shifts, which are quantum objects,  and a measure which is given in terms of the geodesic $X$-ray transform of the potential, which is a classical quantity.  The geodesic $X$-ray transform of a function (or a tensor)  is the map that takes a geodesic to the integral of the function along it.  This extends a result of Bulger and Pushnitski \cite{BulPus1}, proved for Schr\"odinger operators in  Euclidean space, to  hyperbolic space.  

 Such trace formulas have been proved  for  the scattering phase shifts of short range real valued potential perturbations of the Euclidean Laplacian, in different r\'egimes.  First,  Birman and Yafaev \cite{BirYaf1,BirYaf2,Yaf} for fixed energies, while  Bulger and Pushnitski \cite{BulPus1} treated the high energy limit.    The case of magnetic Schr\"odinger operators was considered in \cite{BulPus2} in the high energy limit and by Nakamura \cite{Nak} for fixed energies, but for more general manifolds. We will not treat magnetic potentials here.
 
  In the case studied here as well as in \cite{BulPus1}, the potential is an $O(h^2)$ semiclassical perturbation.  In the case where the potential is part of the semiclassical principal symbol of the Hamiltonian, the distribution of scattering phase shifts for the semiclassical Schr\"odinger operator in Euclidean space was studied by Gell-Redman and Hassell \cite{GelHas0,GelHas},  Datchev, Gell-Redman, Hassell and Humphries \cite{DGR},  and by Gell-Redman, Hassell and Zelditch  \cite{GelHasZel}.

Gell-Redman and Ingremeau \cite{GelIng} studied the problem for scattering by convex obstacles in Euclidean spaces.  Ingremeau \cite{Ing} studied compactly supported metric and potential perturbations, and perhaps more importantly, he replaced the assumption of non-trapping  with a milder condition on the trapped set.    Even in the case of radial potentials, the analysis of the semiclassical problem is very delicate when there is trapping \cite{Ber}. We also mention the work of Christiansen and Uribe \cite{ChrUri} for the semiclassical Schr\"odinger operator in the case of manifolds with one cylindrical end.  
 
  There are several equivalent ways of defining the scattering matrices for the perturbed and unperturbed Laplacians. One can give a dynamical definition in terms of the wave groups and radiation fields as in \cite{Lax,LaxPhi1,LaxPhi2,LaxPhi3,Nak,SaB}, but we will not pursue this here. Instead, we will work in the frequency domain and define the scattering matrix in terms of the asymptotic expansion of generalized eigenfunctions.

  It is important to compare the behavior of the eigenvalues of the scattering matrices in the in Euclidean and hyperbolic spaces.  In Euclidean space,  the unperturbed scattering matrix is the identity, which has only one eigenvalue with infinite multiplicity, independently of the energy.  On the other hand, for short range potentials, and for a fixed energy, the scattering matrix of the perturbed Hamiltonian has a countable set of eigenvalues on $\ms^1$ which possibly accumulate at $(1,0)$, see for example \cite{BirYaf1,BirYaf2,Yaf}.   The scattering matrix of the Laplacian in hyperbolic space $\mhn$, defined in equations \eqref{metg0} and \eqref{defg0} below, has special features.  If $g_0$ is the metric on $\mhn,$  then for each defining function $x$ of the boundary of $\mhn$, $x^2g_0\bigr|_{\{x=0\}}=H$ defines a metric on the boundary $\p \mhn$, which is the sphere at infinity.  A different choice of a defining function of $\p \mhn$,  say $x'= e^{\alpha(x,\theta)} x$, induces a different metric $H'={x'}^2 g_0\bigr|_{x'=0}= e^{2\alpha(0,\theta)} H$. This induces a conformal structure on $\p \mhn$ and the definition of the scattering matrix has to take into account this conformal structure on $\p \mhn$.  For a particular choice of $x$, the scattering matrix is given by equation \eqref{mca0} below and it has a lot of eigenvalues.  One can see that  for a fixed energy the eigenvalues are dense on $\ms^1$, and for variable high energies, they actually cover the entire $\ms^1$.   One defines the relative scattering matrix corresponding to a perturbation to be the product of the scattering matrix corresponding to the perturbation times the inverse of the scattering matrix of the unperturbed Laplacian.  It turns out that, in the case we are dealing with,  the eigenvalues of the relative scattering matrix do not depend on the choice of a conformal representative of the metric on the sphere at infinity.   Moreover, its eigenvalues  only possibly accumulate  at $(1,0)$, which is somewhat the opposite of what happens in the Euclidean case.

 The proof of our main result, Theorem \ref{main}, follows the strategy of \cite{BulPus1}. The main idea is to  use the Born approximation to write the relative scattering matrix in two parts, one of which  is a semiclassical pseudodifferential operator whose principal symbol is equal to the geodesic $X$-ray transform of the potential and then use Schatten-von Neumann estimates to show that the contribution of the second term to the trace formula is of lower semiclassical order.    In order to do that,  we prove Schatten-von Neumann estimates for the Poisson operator in hyperbolic space, which is also known as Helgason Fourier transform  \cite{ItoSat,Hel,Lax,Zel}. These are well known for the Fourier transform in the Euclidean setting and can be found in Chapter 8 of \cite{Yaf} for example, but as far as we can tell, we do not have a reference for them  in the case of hyperbolic space.

 As an application of the trace formula, we show that the high energy limit of the phase shifts uniquely determines a radial  potential which is monotone and super-exponentially decaying with respect to the geodesic distance in $\mhn$.  This result also holds for the Euclidean space using the trace formula proved in \cite{BulPus1}.   Levinson \cite{Lev} has shown that, in Euclidean space,  the phase shifts for all energies determine a radial potential (with additional assumption about its decay), provided its Schr\"odinger operator has no negative eigenvalues.  Bargmann \cite{Bar} has given counter-examples showing that the phase shifts do not determine  radial potentials which have negative eigenvalues. Our result does allow negative potentials whose Schr\"odinger operator may have eigenvalues.

\section{Acknowledgements}

The author is partly supported by the Simons Foundation grant \#848410.  He is very grateful to an anonymous referee for carefully reading the paper, making many suggestions on how to improve the exposition and for insisting that more details should be provided.

\section{The Statement of the Main Results}   

We define  $\mh^{n+1}$, $n\geq 1$, to be  the Riemannian manifold of dimension $n+1$ which  consists of the interior of the ball   
\beq\label{metg0}
\mb^{n+1}=\{w\in \mr^{n+1}: |w|<1\} \text{ equipped with the metric } g_0=\frac{4 dw^2}{(1-|w|^2)^2}.
\eeq 
As in \cite{FefGra,MazMel,Megs} one thinks of $\mhn$ as the interior of a $C^\infty$ manifold with boundary equipped with a metric which is singular at the boundary $\p \mb^{n+1}=\msn$ and has a specific asymptotic expansion in a tubular neighborhood of $\p \mhn$.     The metric $g_0$ defined in \eqref{metg0} is of course singular at $\msn$ and if we use Euclidean polar coordinates $(r,\theta)$, $r=|w|$ and $\theta=\frac{w}{|w|}$, it is given by
 \beq\label{defg0}
 g_0=\frac{4(dr^2+r^2 d\theta^2)}{(1-r^2)^2}.
 \eeq
 If we set  
 \beq\label{g0x}
 \begin{split}
 & x=\frac{1-r}{1+r}, \text{ then } \\
  g_0=  \frac{dx^2}{x^2} & + \frac{(1-x^2)^2}{4} \frac{d\theta^2}{x^2}, \;\ x>0, \\
 &  \text{ and } \bigl(x^2 g_0\bigr)\biggr|_{\{x=0\}}=\frac14 d\theta^2.
\end{split}
 \eeq
The function $x$ is a defining function of the boundary, in the sense that $x\in C^\infty$ near $\p \mhn$   and $x=0$ only at $\p \mhn$ and $dx\not=0$ at $\p \mhn$.  But for example $\wt x=1-r$, or $\wt x= e^{\phi} x$, $\phi \in C^\infty(\mhn)$, would also be defining functions of $\p \mbb$. The induced metric on $\msn$, which is given by $(x^2 g_0)\bigr|_{x=0}$ depends on the choice of $x$. Vice-versa, as noted in \cite{Gra,JosSaB}, the choice of a conformal representative of $x^2 g_0\bigl|_{\p \mhn}$ determines a unique boundary defining function $\vrho$ of $\p \mhn$ such that  $\varrho^2 g_0= d\varrho^2+ H(\varrho,\theta,d\theta)$ near $\{\varrho=0\}$.

Let $\Delta_{g_0}$ denote the Laplacian on $\mhn$. It is known that the spectrum of $\Delta_{g_0}$ is absolutely continuous and equal to 
$[\nsq,\infty)$, see for example \cite{IsoKur,LaxPhi2,McKean,Megs}.   If $V$ is real valued and decaying fast enough at infinity, $\Delta_{g_0}+V$ is self-adjoint and its spectrum consists of an absolutely continuous part $[\nsq,\infty)$ and a finite number of eigenvalues  in $(0, \nsq)$, see for example \cite{JosSaB,MazMel,SaBWan1}.  We shall work with $\Delta_{g_0}-\nsq$ which shifts the continuous spectrum to $[0,\infty)$.

The generalized eigenfunctions of $\Delta_{g_0}-\nsq+V$ are solutions of
\[
(\Delta_{g_0}+V-s(n-s)) u(s,w)=0,
\]
which are bounded but not in $L^2$, and have a prescribed asymptotic behavior at infinity, which in this case is the boundary of $\mhn=\msn$,  
given by powers of the boundary defining function $x$, see  equation \eqref{Poi0} below.  Here  $s=\novt+i\la,$ $\la \in \mr\setminus 0$, and so $s(n-s)=\nsq+\la^2,$ and this choice of spectral parameter is standard in this setting.

 This leads to the definition of  $\mca_V(s)$ and $\mca_0(s)$, which are the scattering matrices corresponding to $V$ and $V=0$, see \eqref{scatg0}. If we pick $x$ as in \eqref{g0x}, we have
\beq
\label{mca0}
\mca_{0}(s)= \frac{\Gamma(-i\la)}{\Gamma(i\la)} \Delta_{\ms^n}^{i\la}, \;\ s=\novt+i\la,  \la\in \mr\setminus 0,
\eeq
where $\Delta_{\msn}$ is the standard Laplacian on $\msn$, see for example \cite{GuiZwo,JosSaB,Per}.

The operators $\mca_0(s)$ and $\mca_V(s)$, which in principle act on $C^\infty(\msn)$, extend to unitary operators acting on $L^2(\msn)$ with respect to the metric induced by $g_0$ at the boundary, and hence their eigenvalues lie on $\ms^1$, but depend on the choice of $x$.   To overcome this difficulty, we work with the relative scattering matrix $\Tau_V(s)= \mca_V(s)\mca_0(s)^{-1}$. It turns out that  $\Tau_V(s)$ is unitary and its eigenvalues do not depend on the choice of $x$. Moreover, $\Tau_V(s)= I + \mce_V(s)$, where $\mce_V(s)$ is a compact operator. Therefore, its eigenvalues lie on the unit circle and only possibly accumulate at $(1,0)$.

Since the eigenvalues of the Euclidean Laplacian on $\msn$ are 
\[
 k(k+n-1) \text{ with multiplicity } \binom{n+k}{n},
 \]
 the  spectrum of $\mca_0(\la)$  given by \eqref{mca0} consists of eigenvalues
 \beq\label{eigenA0}
 \mu_k(\la)=\frac{\Gamma(-i\la)}{\Gamma(i\la)} \left( k(k+n-1)\right)^{i\la} \text{ with multiplicity } \binom{n+k}{n}.
 \eeq
 
 For $\la$ fixed,  modulo the term involving the gamma function which is independent of $k$, these are points on $\ms^1$  of the form $e^{i\la \log\left( k(k+n-1)\right)}$ and it is well known, see Example 2.4 of \cite{KuiNie},  that this is  a dense subset of $\ms^1$, but which is not equidistributed.    On the other hand for fixed $k$,  Stirling's formula \cite{Handbook} gives  that
 \[
 \frac{\Gamma(i\la)}{\Gamma(-i\la)}=i\left(\frac{\la^2}{e^2}\right)^{i\la}(1+ O(\frac{1}{\la}))= i e^{2i\la(\log \la-1)}(1+O(\frac{1}{\la})), \text{ as } \la\rightarrow \infty,
 \]
 and therefore as $\la$ varies, it covers the entire $\ms^1$.   If one just takes $\la \in \mn$, then according to Example 2.8 of \cite{KuiNie},  the sequence 
 \[
 \mu_k(\la) = i e^{i(2\la(\log \la-1)+\la \log(k(k+n-1)))}(1+O(\frac{1}{\la})),
 \]
 is equidistributed on $\ms^1$. 
 
Therefore, while for $\la$ fixed the eigenvalues of $\mca_0(s)$ are dense on $\ms^1$,   the eigenvalues of $\Tau_{V}(s)$ only possibly accumulate at the point $(1,0)\in \ms^1$.  At the high energy limit, the eigenvalues of $\mca_0(s)$ cover $\ms^1$  and our goal is to describe the effect potential perturbations have on these eigenvalues.  We denote the eigenvalues of $\Tau_V(s)$ by
\[
\begin{split}
 e^{i \delta_j(s)}, \;\ \del_j(s)\in [-\pi,\pi), \;\ j\in \mn, &  \text{ counted with multiplicity, and we know that } \\
& \del_j(s)\rightarrow 0  \text{ as } j\rightarrow \infty.
\end{split}
\]
The numbers $\del_j(s)$ will be called the {\it relative scattering phase shifts}.    We will show that 
\[
||\Tau_V(s)-\id||_{\mcl(L^2(\mhn))}\leq C|\la|^{-1} \text{ for } s=\novt+i\la, \;\ \la \in \mr, \;\  |\la|>C,
\]
see equation  \eqref{est-EV1} below. So it follows that  $ |e^{i\del_j(s)}-1| \leq C |\la|^{-1},$ for  $\la \in \mr$ and $|\la|>C$. In particular, for $\la\in \mr\setminus 0$, $|\la|>>1$,   $|\del_j(\novt+i\la)|\leq |\sin\del_j(\novt+i\la)|\leq C |\la|^{-1}$.  We set $\la=\ooh$ and for $s=\novt+\frac{i}{h}$,  we denote  $\ka_j(h)=\del_j(\novt+\ioh)$.  As in \cite{BulPus1}, we define the following measure
\beq\label{measNew}
\lan \mu_h,f\ran= (2\pi h)^{n} \sum_{e^{i\ka_j(h)}\in \operatorname{spec} \Tau_V} f\left(\frac{\ka_j(h)}{h}\right), \;\ f\in C_0^0(\mr \setminus 0).
\eeq

If one takes $f=\chi_{([a,b])}$, $0\not \in [a,b]$, (as a limit of continuous functions),  this is equivalent to the following counting measure:
\beq\label{count2}
(2\pi h)^{-n} \lan \mu_h,\chi_{[a,b]}\ran= \#\{j: \frac{\ka_j(h)}{h} \in [a, b]\subset \mr \setminus 0\},  \text{ counted with multiplicity.}
\eeq

Our purpose is to study the limit of  $\mu_h$, as $h\downarrow 0$ and we prove the following
\begin{theorem}\label{main}   Let $m\in (1,\infty)$ and $n\in \mn$, $n\geq 1$. Let $V\in x^mL^\infty(\mhn)$  be real valued, and let 
$f \in C^0(\mr)$ be such that  $\frac{f(t)}{t^p} \in C^0(\mr)$, for some  $p\in \mn$  with $p> \frac{2(n-\ha)}{m-1}$.
 Then
\beq\label{limh1}
 \lim_{h\rightarrow 0}  \lan \mu_h,f\ran = \int_{\ms^n} \int_{\mrn} f\biggl( -2^{n-1} X(V)\bigl(\frac{\xi}{2},\theta\bigr)\biggr) \ d\xi d\theta. 
\eeq
  If $V$ is bounded and compactly supported,  we can take $m=\infty$ and $p=1$.  Equation \eqref{limh1}  can be rephrased as 
\beq\label{lim-high2}
 \lim_{h\rightarrow 0}  \lan \mu_h,f\ran = \int_{\mr} f(t) \ d\nu(t), \;\ f\in C_0^\infty(\mr\setminus 0),
\eeq
where the measure $\nu$ is defined to be
\beq\label{meas3}
\begin{split}
 \nu(\alpha,\beta)= & \operatorname{meas}\{(\theta,\xi), \; \theta\in \msn, \xi\in \mrn, \; \lan \xi, \theta\ran=0 \text{ such that } -2^{n-1} X(V)(\frac{\xi}{2},\theta)\in (\alpha,\beta) \subset \mr\setminus 0\}.
\end{split}
\eeq
Here  ``meas'' indicates the Lebesgue measure, and $ X(V)(\frac{\xi}{2},\theta)$ is the geodesic $X$-ray transform of the potential $V$ along the geodesic $\sigma(\frac{\xi}{2},\theta)$ defined in \eqref{RTV} below.
\end{theorem}
As already mentioned, in the case of a potential perturbation of $\Delta$ in $\mrn$, the analogue of this theorem was proved  in \cite{BulPus1}.

\subsection{An Inverse Theorem}

We apply Theorem \ref{main} to prove that the high energy limit of the phase shifts determines a class of radial potentials in $\mhn$.  The methods of \cite{Lev} will likely also work in $\mhn$,  but we are not aware this has been investigated.  However, we are only dealing with the high energy limit of the phase shifts and we only recover the measure $\nu$, so it is reasonable to expect that we need to restrict the class of potentials.  

We say that a potential $V(w)$, $w\in \mhn$, is radial if it only depends on the distance $d_{g_0}(0,w)$ from the origin to $w$ with respect to the metric $g_0$.  It follows form equation \eqref{dist0} below that $\rho(w)=d_{g_0}(w,0)=\log \frac{1+|w|}{1-|w|}$, where $|w|$ denotes the Euclidean norm of $w$.  So $|w|=\tanh(\frac{\rho}{2})$. So if $V$ is a radial potential if it can be written it as
\[
V(w)= V(\rho(w))=V\biggl( \log \frac{1+r}{1-r}\biggr), \;\ r=|w|
\]
Moreover, if  $w=\theta+ \frac{t\theta+ \xi}{t^2+|\xi|^2}$, then $|w|^2=1+ \frac{1+2t}{t^2+|\xi|^2}$ (see the discussion in Section \ref{GRT} to understand the change of variables), and so it follows from \eqref{RTV} that $X(V)(\xi,\theta)$ is radial and we define
\beq\label{defgv}
-2^{n-1} X(V)(\frac{\xi}{2},\theta)\defn G_V(r), \;\ r=|\xi|.
\eeq

\begin{theorem}\label{inverse}  Let $V_1, V_2 \in C^1(\mhn)\bigcap x^m L^\infty(\mhn)$, for all $m>0$, be real valued radial potentials which are strictly monotone in the sense that $V_j'(r)\not=0$ if $V_j(r)\not=0$, $j=1,2$.  Let  $\mu_h(V_j)$, $j=1,2$ be the measure defined in \eqref{measNew}.  If
\beq\label{limit-meas}
\lim_{h\rightarrow 0} \lan \mu_h(V_1), f\ran= \lim_{h\rightarrow 0} \lan \mu_h(V_2), f\ran, \text{ for all } f\in C_0^\infty(\mr\setminus 0),
\eeq
then $V_1(w)=V_2(w)$ for all $w\in \mhn$.
\end{theorem}

Theorem \ref{inverse} allows negative potentials, and the associated Schr\"odinger operators  may in principle have discrete spectrum, and so it covers cases not included in the (possible generalization of the) results of \cite{Lev}.
The proof of Theorem \ref{inverse} also applies to the Euclidean case and it follows from the trace formula proved in \cite{BulPus1}. 

 It follows from  \eqref{limit-meas}  that $\nu_{{}_{V_1}}=\nu_{{}_{V_2}}$, where $\nu_{{}_{V_j}}$ is the measure defined by \eqref{meas3} corresponding to the potential $V_j$. Notice that these are invariant under measure preserving diffeomorphisms of $\msn \times \mrn$ which also preserve orthogonality.  In the case of radial potentials, these become measures on $\mr$ and  the best one can hope to conclude is that $G_{V_1}(r)=G_{V_2}(r+C)$. We prove that $G_{V_1}$ and $G_{V_2}$ have the same range and since they are monotone, it follows that $G_{V_1}(0)=G_{V_2}(0),$ but then it follows that $G_{V_1}(0)=G_{V_2}(0)= G_{V_2}(C)$ and this implies that $C=0.$
 
 But since  $G_{V_1}=G_{V_2}$  and the potentials decay super-exponentially,  we can conclude that $V_1=V_2$, because the $X$-ray transform is injective for functions in $x^m L^\infty(\mhn)$ for all $m>0$, see Theorem 1.2 and Corollary 1.3  of  Chapter 3 of \cite{Hel}.

\section{An Outline of the Proof of Theorem \ref{main}}\label{OUT}

First of all, one needs to recognize that when $f(t)=t^k$, $k\in \mn$, Theorem \ref{main} is a result about the rescaled trace of the operator $(\Tau_V(\novt+\ioh)-\id)^k$.  Formally speaking, we first notice that 
\[
(2\pi h)^n Tr \biggl( \frac{1}{h}\im \biggl(\Tau_V\bigl(\novt+\ioh\bigr)-\id \biggr)\biggr)^k= (2\pi h)^n \sum_{e^{i\ka_j(h)}\in \spec \Tau_V} \biggl(\frac{\sin \ka_j(h)}{h}\biggr)^k, \;\ \ka_j(h)=\del_j(\novt+\ioh),
\]
 but we use methods of \cite{BulPus1} to show that if we replace $\sin \ka_j(h)$ with $\ka_j(h)$ in the formula, the difference between the two is a term of order $O(h)$. We also show that the imaginary part of the operator can be replaced by  the operator itself, again with an $O(h)$ error.  So when $f(t)=t^k$,  the trace formula says that
\[
(2\pi h)^n Tr \biggl( \frac{1}{h}\biggl( \Tau_V\bigl(\novt+\ioh\bigr)-\id \biggr)\biggr)^k= \int_{\msn}\int_{\mrn} \biggl(-2^{n-1} X(V)(\frac{\xi}{2}, \theta)\biggr)^k \ d\xi d\theta+ O(h).
\]
For more general functions, the result follows from an application of the Stone–Weierstrass theorem.  The idea of proving a trace formula for polynomials first and then using the Stone–Weierstrass theorem to extend it to a more general class of functions is certainly not new and it has been used in several contexts. In the particular case of scattering shifts  in appears in several papers, see for example \cite{BirYaf1,BirYaf2,BulPus1,DGR,GelHasZel,GelIng,Ing}.

The trace of the operator  $\Tau_V(\novt+\ioh)-\id$  is equal to the integral of the its Schwartz kernel restricted to the diagonal, and our proof relies on the fact that the Schwartz kernels of the scattering matrices $\mca_0(s)$ and $\mca_V(s)$ can be obtained from the kernels of the respective resolvents \cite{GuiZwo,JosSaB}.  In this paper, as in \cite{BulPus1},  we use the Born approximation for $\Tau_V(\novt+\ioh)$ to show that  $\Tau_V(\novt+\ioh)-\id= \mcu_V(\novt+\ioh)-\mcw_V(\novt+\ioh)$, and that for bounded compactly supported $V$, $\frac{1}{h}\mcu_V(\novt+\ioh)$ is a semiclassical pseudodifferential operator with principal symbol  $-2^{n-1} X_V(\frac{\xi}{2},\theta)$. We do not prove that  $\mcw_V(\novt+\ioh)$ is a pseudodifferential operator, and instead we use Schatten--von Neumann norm estimates to show that its contribution to the trace is of order $O(h)$, namely:
\[
(2\pi h)^n Tr \biggl( \frac{1}{h}\biggl( \Tau_V\bigl(\novt+\ioh\bigr)-\id \biggr)\biggr)^k- (2\pi h)^n Tr\biggl( \frac{1}{h} \mcu_V\bigl(\novt+\ioh\bigr)\biggr)^k= O(h).
\]
 The Scatten-von Neumann estimates are also needed to extend the result to non-compactly supported short range potentials.

 As we will see in Section \ref{SBP} below, a parametrix for the resolvent gives a parametrix for the scattering matrix, and one should expect to be able to use a more precise parametrix for the semiclassical resolvent such as the ones constructed in \cite{CheHas,MelSaBVas,SaBWan} to directly show  that if $V$ is compactly supported, $\ooh(\Tau_V-\id)$ is a semiclassical pseudodifferential operator with principal symbol  $-2^{n-1} X_V(\frac{\xi}{2},\theta)$.
  Then Theorem \ref{main}  would essentially follow from the calculus of semiclassical pseudodifferential operators.   This approach would also possibly work  in more general situations, including non-trapping asymptotically hyperbolic manifolds, but it would be much more refined and hopefully will be pursued elsewhere.

\section{The  Geodesic $X$-Ray Transform in $\mhn$}\label{GRT}

It is well known, see for example \cite{Hel,Laf}, that the geodesics of $\mhn$ are either circles which intersect the boundary $\msn$ perpendicularly, or straight lines that go to through the origin.  Of course there are many different ways of describing such curves, but in this paper they will appear naturally in the computation of the principal symbol of the first Born approximation of the semiclassical scattering matrix by using the inversion map with respect to the sphere $\msn:$
\[
\begin{split}
 Q:   \mr^{n+1} & \setminus 0 \longrightarrow \mr^{n+1}\setminus 0, \\
& w \longmapsto \frac{w}{|w|^2}.
\end{split}
\]

\begin{figure}[h!]
\begin{tikzpicture}[scale=.5]
\centering
\draw (0,0) circle (71pt);
\draw (0,-4) -- (0,4);
\draw (-4,0) -- (4,0);
\draw[thick] (0,0) -- (2.5,0) node[right,below]{$\;\ \;\ \theta$};
\draw[thick] (2.5,0) -- (4,0) node[right,below]{$\;\ \;\ t \theta$};
\draw[thick] (4,0) -- (5,0); 
\draw[thick] (0,0) -- (0,3.2) node[blue,right]{$ \xi_1\;\ \;\ $};
\draw[thick] (0,0) -- (0,2.2) node[red,right]{$ \xi_2\;\ $};
\draw[thick] (4,0) -- (4,3)  node[sloped,right]{\;\ \;\ $t\theta+\blue{\xi_1}$};
\draw[thick] (4,0) -- (4,2)  node[sloped,right]{\;\ \;\ $t\theta+\red{\xi_2}$};
\draw[blue] (0,0.5) circle(14pt);
\draw[red] (0,3) circle(85pt);
\draw[blue] (-4,3)--(5,3);
\draw[red] (-4,2)--(5,2);
\draw (0,0)--(4,3);
\draw[thick] (0,0)--(0.5,0.4);
\tkzDefPoint(2.5,0){A};
\tkzDefPoint(0.5,0.4){B};
\tkzDefPoint(4,0){C};
\tkzDefPoint(4,3){D};
\tkzDefPoint(0,3){E};
\tkzDefPoint(0,2){F};
\tkzDefPoint(4,2){G};
\foreach \n in {A,B,C,D,E,F,G}
  \node at (\n)[circle,fill,inner sep=1.5pt]{};
\end{tikzpicture}
\caption{Fixed $\theta\in \ms^1$, the blue (smaller) and red (larger) circles are the images of the lines $t\longmapsto t\theta+\xi_j$,  $\xi_j\in \mr^{2}\setminus 0$, $ \lan \xi_j,\theta\ran=0$, $j=1,2$, respectively, with respect to the map $Q$. The circles are oriented clockwise and as $\pm t\rightarrow \infty$ the point $Q(t\theta+\xi)$ approaches the origin from the right and from the left, respectively.}
\label{circ1}
\end{figure}
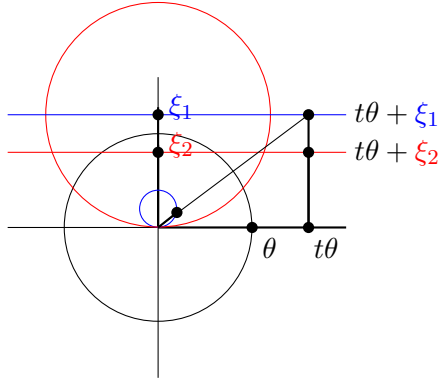

One can use this map to parametrize all such circles. Let  $\theta \in \msn$ and $\xi\in \mr^{n+1}\setminus 0$, with $\lan \xi, \theta\ran=0$, and define  the curve
\[
t \longmapsto  Q(t\theta+\xi)= \frac{t\theta+\xi}{t^2+|\xi|^2}, \text{ where } \lan \xi,\theta\ran=0,
\]
obtained by inverting the line $t\mapsto t\theta+\xi$ with respect to $\msn$.
Notice that
\[
\left| \frac{t\theta+\xi}{t^2+|\xi|^2}-\frac{\xi}{2|\xi|^2}\right|^2=\frac{1}{4 |\xi|^2}, \text{ provided }  \xi\not=0, \text{ and } \lan \xi,\theta\ran=0,
\]
and so if we think of $\theta=(1,0,\ldots,0)$ and $\xi=(0,\xi_1,\ldots, \xi_{n+1})$,  this is a parametrization of the circle centered at $(0,\frac{\xi}{2|\xi|^2})$ and radius $\frac{1}{2|\xi|};$ the case $\xi=0$ represents a ray in the $\theta$ direction.  The closure of these circles go through the origin at the limits $t\rightarrow \pm \infty$ and  they are tangent to the $\theta$-axis and intersect the $\xi$-axis orthogonally.  Moreover, since they are contained on the hyperplane defined by $\xi,\theta$ passing through the origin, we we just need to visualize them in two dimensions, see Fig.\ref{circ1}.

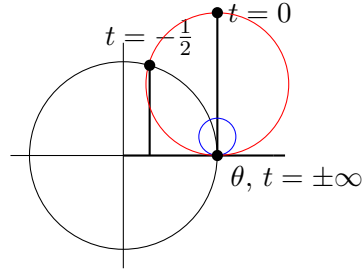
\begin{figure}[h!]
\begin{tikzpicture}[scale=.5]
\centering
\draw (0,0) circle (71pt);
\draw (0,-3) -- (0,3);
\draw (-3,0) -- (3,0);
\draw[thick] (0,0) -- (2.5,0) node[right,below]{\;\ \;\ \;\ \;\ \;\ \;\ \;\ \;\ \;\ $\theta$, $t=\pm\infty$};
\draw[thick] (0.7,0) -- (0.7,2.4) node[left,above]{$t=-\ha$};
\draw[thick] (2.5,0) -- (2.5,3.8) node[right]{$t=0$};
\draw[thick] (2.5,0) -- (4,0);
\draw[thick] (4,0) -- (4.3,0);
\draw[blue] (2.5,0.5) circle(14pt);
\draw[red] (2.5,1.9) circle(54pt);
\tkzDefPoint(2.5,0){A};
\tkzDefPoint(2.5,3.8){B};
\tkzDefPoint(0.7,2.4){C};
\foreach \n in {A,B,C}
 \node at (\n)[circle,fill,inner sep=1.5pt]{};
\end{tikzpicture}
\caption{The circles  of Fig.\ref{circ1} shifted by $\theta$.}
\label{circ2}
\end{figure}

 Of course, the map 
\beq\label{great-circ}
t \longmapsto \ga(t,\xi,\theta)=\theta+ \frac{t\theta+\xi}{t^2+|\xi|^2},
\eeq
just shifts the circle by $\theta$ and so it is tangent to the $\theta$ axis and  centered at $(\theta, \frac{\xi}{2|\xi|^2})$.  Since $\lan \xi, \theta\ran =0$,  the circle \eqref{great-circ} intersects $\msn$ orthogonally at $\theta$, see Fig.\ref{circ2}.

Notice that 
\[
\left| \theta+ \frac{t\theta+\xi}{t^2+|\xi|^2}\right|^2= 1+ \frac{1+2t}{t^2+|\xi|^2},
\] 
and so these circles intersect $\msn$ when $t=-\ha$. Moreover, the tangent vector
\[
T(t,\xi,\theta)=\p_t \ga(t,\xi,\theta)=\frac{(|\xi|^2-t^2) \theta -2t\xi}{(t^2+|\xi|^2)^2},
\]
satisfies
\[
T(-\ha,\xi,\theta)= \frac{4}{1+|\xi|^2} \ga(-\ha,\xi,\theta),
\]
and so $\ga(t,\xi,\theta)$ also intersects $\msn$ orthogonally when $t=-\ha$.  Therefore, for $t\in (-\infty,-\ha)$, $\ga(t,\xi,\theta)$ is a geodesic of $\mhn$, and all geodesics on $\mhn$ can be parametrized in this way, and we shall denote them by $\sigma(\theta,\xi)$. 

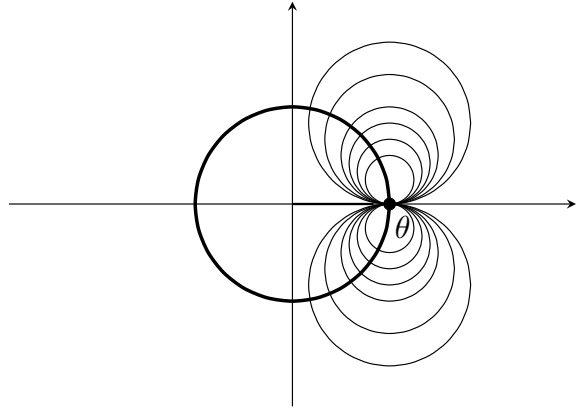
\begin{figure}
\begin{tikzpicture}[x=1mm,y=1mm,scale=1.1]
\centering
    \begin{axis}[
            xmin=-3.5,xmax=3.5,ymin=-2.5,ymax=2.5,
            axis lines=center,
            ticks=none,
            unit vector ratio*=1 1 1,
           ]        
          \addplot [smooth, domain=(2*pi:0)] ({1.2+0.8*cos(deg(x))},{0.8+ 0.8*sin(deg(x))});         
           \addplot [smooth, domain=(2*pi:0)] ({1.2+0.8*cos(deg(x))},{-0.8+ 0.8*sin(deg(x))});         
          \addplot [smooth, domain=(2*pi:0)] ({1.2+0.3*cos(deg(x))},{0.3+ 0.3*sin(deg(x))});         
           \addplot [smooth, domain=(2*pi:0)] ({1.2+ 0.3*cos(deg(x))},{-0.3+ 0.3*sin(deg(x))});         
           \addplot [smooth, domain=(2*pi:0)] ({1.2+0.4*cos(deg(x))},{0.4+ 0.4*sin(deg(x))});         
            \addplot [smooth, domain=(2*pi:0)] ({1.2+0.4*cos(deg(x))},{-0.4+ 0.4*sin(deg(x))});         
             \addplot [smooth, domain=(2*pi:0)] ({1.2+0.5*cos(deg(x))},{0.5+ 0.5*sin(deg(x))});         
             \addplot [smooth, domain=(2*pi:0)] ({1.2+0.5*cos(deg(x))},{-0.5+ 0.5*sin(deg(x))});         
              \addplot [smooth, domain=(2*pi:0)] ({1.2+0.6*cos(deg(x))},{-0.6+ 0.6*sin(deg(x))});         
               \addplot [smooth, domain=(2*pi:0)] ({1.2+0.6*cos(deg(x))},{0.6+ 0.6*sin(deg(x))});         
                \addplot [smooth, domain=(2*pi:0)] ({1.2+cos(deg(x))},{1+sin(deg(x))});         
                 \addplot [smooth, domain=(2*pi:0)] ({1.2+cos(deg(x))},{-1+ sin(deg(x))});         
          \addplot [very thick, smooth, domain=(0:2*pi)] ({1.2*cos(deg(x))},{1.2*sin(deg(x))});         
 \addplot[thick,smooth,domain=(0:1)]({x},{0})node[right,below]{\;\ \;\ \;\ $\theta$};
        \addplot[only marks,mark=*] coordinates{(1.2,0)};
         \end{axis}
\end{tikzpicture}
\caption{A family of circles $\ga(t,\xi,\theta)$, for  $\theta$ fixed and for several values of $\xi\not=0$, with $\lan \xi,\theta\ran=0$, as $\xi$ points upward or downward.}
\label{circ3}
\end{figure}

This defines a map of the tangent bundle of $\msn$ to itself
\[
\begin{split}
\mcs: T\msn \longrightarrow T\msn, \\
(\theta,\xi) \longmapsto (\theta',\xi'),
\end{split}
\]
such that $\ga(-\ha,\xi,\theta)=\theta'$ and $\ga(-\ha, -\xi',\theta')=\theta$. In other words, one follows the geodesic $\ga(t,\xi,\theta)$, which starts at $\theta$ in the direction determined by  $\xi$ until one intersects $\msn$ again when $t=-\ha$, at  $\theta'$. The vector $\xi'$ is defined to be such the geodesic starting at $\theta'$ defined by $-\xi'$ going backwards  intersects $\msn$ at $\theta$.  One can define this map, or perhaps more appropriately,  its dual on the cotangent bundle of $\msn$, as the scattering relation on $\msn$.

The length of the tangent vector with respect to the Euclidean metric is
\[
 \left| \p_t\ga(t,\xi,\theta)\right|=\frac{1}{t^2+|\xi|^2}
\]
and so if $t\in (-\infty,-\ha)$, then $|\ga(t,\xi,\theta)|<1$, and  the length of the tangent vector with respect to the metric $g_0$ at the point $\theta+ \frac{t\theta+\xi}{t^2+|\xi|^2}$ is given by
\[
 \left| \p_t\ga(t,\xi,\theta)\right|_{g_0}=2 \left(1-\left|\theta+\frac{t\theta+\xi}{t^2+|\xi|^2}\right|^2\right)^{-1}|\p_t \ga(t,\xi,\theta)|= \frac{-2}{1+2t}, \;\ t\in (-\infty,-\ha).
 \]

Notice that this is positive for $t\in (-\infty,-\ha)$, and blows up as $t\rightarrow -\ha$, where it intersects the boundary, since the metric $g_0$ also blows up there.
  
   The geodesic $X$-ray transform of a function $U\in C_0^\infty(\mhn)$ is defined to be the map which takes a geodesic $\ga$ of $\mhn$ into the integral of $U$ along $\ga$, see for example Chapters 1 and 3 of \cite{Hel}.  We parametrize such geodesics as in \eqref{great-circ} and so,  we  define the geodesic $X$-ray transform of $U$ along the geodesic $\ga(\theta,\xi)$ to be
\beq\label{RTV}
\begin{split}
& X(U)(\xi,\theta)=  \int_{\ga(\theta,\xi)} U(w) ds= 
 \int_{-\infty}^{-\ha} U\bigl(\theta+ \frac{t\theta+\xi}{t^2+|\xi|^2}\bigr) \left| \p_t\left(\theta+ \frac{t\theta+\xi}{t^2+|\xi|^2}\right)\right|_{g_0} \ dt= \\
&- \int_{-\infty}^{-\ha} U\bigl(\theta+ \frac{t\theta+\xi}{t^2+|\xi|^2}\bigr) \frac{2}{1+2t} \ dt, \;\ \theta\in \msn, \;\ \xi\in \mr^{n+1} \text{ such that } \lan \xi,\theta\ran=0.
 \end{split}
 \eeq
 Notice that if $U(w)$ is compactly supported in $\mhn$, so is $X(U)(\xi,\theta)$.  To see that, suppose that $U(w)=0$ if $|w|^2>1-\del$. This implies that $U\biggl(\theta+\frac{t\theta+\xi}{t^2+|\xi|^2}\biggr)=0$ if $
 \bigl|\theta+\frac{t\theta+\xi}{t^2+|\xi|^2}\bigr|^2-1= \frac{2t+1}{ t^2+|\xi|^2} > -\del$. But this is equivalent to
 \[
 \bigl(t+\frac{1}{\del}\bigr)^2+ |\xi|^2+ \frac{1}{\del}-\frac{1}{\del^2}>0.
 \]
It follows that $U(\theta+ \frac{t\theta+\xi}{t^2+|\xi|^2}) =0$ if $|\xi|\geq \frac{1}{\del}$ and in particular, $X(U)(\xi,\theta)=0$ if $|\xi|\geq\frac{1}{\del}$.
 
 But $X(U)$ is also well defined for any $U$ which is equal to zero outside the ball $\mb^{n+1}$ and for which the integral converges. For instance, we will work with the class of functions  $U\in x^m L^\infty(\mhn)$, $m>1,$  and we denote 
 \beq\label{norm1}
 ||U||_{\infty,m}=\sup_{w\in \mbb} | x(w)^{-m}U(w)|.
 \eeq
 These norms are equivalent for different choices of the boundary defining function $x$ and
 \beq\label{bdXR}
 |X(U)(\xi,\theta)| \leq C ||U||_{\infty,m} \int_{-\infty}^{-\ha} \frac{|1+2t|^{m-1}}{(t^2+|\xi|^2)^m}\ dt\leq C_m ||U||_{\infty,m} (1+|\xi|)^{-m}.
 \eeq

\section{Scattering by a Potential  in $\mhn$}\label{SBP}

 We discuss some basic facts about scattering theory on hyperbolic spaces, and we refer the reader to \cite{IsoKur,GuiZwo,Lax,LaxPhi1,LaxPhi2,LaxPhi3,Megs,MazMel,Per} for a more thorough discussion and generalization to asymptotically hyperbolic manifolds.  As already mentioned, throughout the paper the potential $V$ is assumed to be real valued.

Let $\Delta_{g_0}$ denote the (positive) Laplacian with respect to the metric $g_0$. Its spectrum is absolutely continuous and equal to $[\nsq,\infty)$,  see for example \cite{IsoKur,LaxPhi2,McKean}.  If $V\in x^m L^\infty(\mhn)$, $m\in [1,\infty)$, the spectrum of  $\Delta_{g_0}+V$,  has two parts: a point spectrum consisting of finitely many eigenvalues in $(0,\nsq)$ and absolutely continuous part which  is also equal to $[\nsq,\infty)$. There are no embedded eigenvalues, including the bottom of the continuous spectrum, see for example \cite{Bou,JosSaB,MazMel,Maz2}.

We consider the generalized eigenfunctions of $\Delta_{g_0}+V-\nsq$. Here, as standard, we use $s=\novt+i\la$, and so $s(n-s)=\la^2+\frac{n^2}{4}$.  It is known, see for example \cite{MazMel,Megs,JosSaB},  that for $V\in x^m L^\infty(\mhn)$, $m\in[1,\infty)$, and any $f\in C^\infty(\ms^n)$ and $\la \in \mr\setminus 0$,  there exists a unique $u(s,w) \in C^\infty(\mhn)$ such that for $x$ as defined in \eqref{g0x}, 
\beq\label{Poi0}
\begin{split}
& \left(\Delta_{g_0}+V-s(n-s))\right)u(s,w)=0 \text{ in } \mhn \text{ and } \\
u(s,w)= \, &  x^{s} f(\theta) + x^{n-s} g(s,\theta)+ O(x^{\novt+1}), \;\ w=(x,\theta) \text{ as } x\searrow 0.
\end{split}
\eeq 

The map
\beq\label{Poi1}
\begin{split}
\mcp_{V}(s): \, & C^\infty(\ms^n) \longrightarrow C^\infty(\mhn) \\
& \mcp_{V}(s) f = u(s,w),
\end{split}
\eeq
for obvious reasons, is called the Poisson operator, and the map
\beq\label{scatg0}
\begin{split}
\mca_{V}(s): \, & C^\infty(\ms^n) \longrightarrow C^\infty(\ms^n) \\
& \mca_{V}(s) f=g,
\end{split}
\eeq
is called the Scattering matrix.  We shall use $\mcp_0(s)$ and $\mca_0(s)$ to denote the Poisson operator and the scattering matrix for $V=0$. 

For $\la\in \mr\setminus 0$, $s=\novt+i\la$,  it follows from the definition that
\beq\label{minusL0}
 \begin{split}
 \mcp_{V}(s)&  \mca_{V}(n-s)= \mcp_{V}(n-s) \text{ and } \\
&  \mca_V^{-1}(s)= \mca_V^{*}(s)= \mca_V(n-s).
\end{split}
\eeq 
and so if  $s=\novt+i\la$, $\la\in \mr\setminus 0$, $\mca_V(s)$ extends to a unitary operator on $L^2(\msn)$ with respect to the metric induced on $\msn$ by the choice of the boundary defining function $x$, as in \eqref{g0x}.   

  It was shown in \cite{JosSaB},  in greater generality,  that for fixed $\la\in \mr\setminus 0$,  $s=\novt+i\la$,  and the same choice of $x$, $\mca_0(s)$ and $\mca_V(s)$ are pseudodifferential operators of order zero (truly of complex order $2i\la$), and moreover
\beq\label{prin-part}
\begin{split}
& \mca_V(s)= \mca_0(s)(I + \mce_V(s)),  \text{ where } \mce_V(s) \in \Psi^{-1}(\msn).
\end{split}
\eeq
This notation means that  $\mce_V(s)$ is a pseudo-differential operator of order $-1$ on $\msn$.   This shows that $\mca_V(s)$ and $\mca_0(s)$ have the same principal part and therefore
\beq\label{ps-rel-sm}
\Tau_V(s)=\mca_V(s) \mca_0(s)^{-1}= I+ \mce_V(s), \;\  \mce_V(s) \in \Psi^{-1}(\msn), \;\ s=\novt+i\la, \;\ \la \in \mr\setminus 0.
\eeq

Since $\mce_V(s)$ is a pseudodifferential operator of order $-1$, and $\msn$ is compact,  then by the Sobolev embedding theorem
\[
\mce_V(s): L^2(\msn) \longmapsto H^1(\msn) \hookrightarrow L^2(\msn)
\]
is a compact operator  and therefore the spectrum of $\Tau_V(s) =\id + \mce_V(s)$  consists of eigenvalues contained in $\ms^1$ which possibly accumulate  at $(1,0)$.  We shall define $\Tau_V(s)$ to be the {\it relative scattering matrix.} This particular relative scattering matrix was also considered by Borthwick and Wang \cite{BorWan} to study the existence of resonances of $\Delta_{g_0}+V$.

 Notice that, if we choose two different boundary defining functions $x$ and  $\wt x$ such that $x= e^{\vphi(\wt x,\theta)} \wt x$, then it follows from \eqref{Poi0} that  the corresponding scattering matrices satisfy
\[
 \mca_\bullet(s)f = e^{(s-n)\vphi(0,\theta)}\wt \mca_\bullet(s) \left(e^{s\vphi(0,\theta)} f\right), \;\ \bullet=V, 0.
 \]
 Since $\mca_0^{-1}(s)=\mca_0(n-s)$, it follows that
\[
 \Tau_V(s)=  e^{(s-n)\vphi(0,\theta)}\wt \mca_V(s) \wt\mca_0^{-1}(s) e^{(n-s)\vphi(0,\theta)}= e^{(s-n)\vphi(0,\theta)}\wt \Tau_V(s) e^{(n-s)\vphi(0,\theta)},
\]
and  so $\Tau_V(s)$ and $\wt \Tau_V(s)$ have the same eigenvalues with the same multiplicity, and we conclude that the eigenvalues of $\Tau_V(s)$ are independent of the choice of $x$.  

\subsection{The Kernels of the Poisson Operator and the Scattering Matrix}   It is a standard convention to define the resolvent in hyperbolic space as 
 \[
 R_{V}(s)=\bigl(\Delta_{g_0}+V- s(n-s)\bigr)^{-1}.
 \]
where $s=\novt+i\la$ so $\nsq+\la^2=s(n-s)$.  When $V$ decays fast enough,  $R_V(s)$ is a holomorphic family of bounded operators  for $\im \la<<0$, or $\re s>> \novt$. In fact, if $V\in x^m L^\infty(\mhn)$, $m\in [1,\infty)$,  $R_V(s)$ continues meromorphically to $\mathbb{C}$, see for example \cite{JosSaB,MazMel}, with poles of finite multiplicity.  These poles are called resonances.  In the case where $V\in C_0^\infty(\mhn)$,  Borthwick \cite{Bor1} and Borthwick and Crompton \cite{BorCro},  proved upper bounds for the counting function of the resonances. Borthwick and Wang \cite{BorWan} proved the existence of resonances and obtained lower bounds for their counting function.

 The Schwartz kernel of the resolvent $R_{\bullet}(s)$, $\bullet=0,V$,  is a distribution in $ \mcd'(\mb^{n+1} \times \mb^{n+1})$,   which we denote  by $R_{\bullet}(s, w,w')$. We know, see for example Lemma 2.1 of  \cite{GuiZwo}, that in the case $V=0$, the Schwartz kernel $R_{0}(s,w,w')$  is given by
\beq\label{ker-R0}
\begin{split}
R_0(s,w,w')&= ( \cosh (d_{g_0}(w,w')))^{-s} G(s,\cosh (d_{g_0}(w,w'))), \;\ \text{ and for } \tau>1, \\
& G(s,\tau)= \pi^{-\novt} 2^{-s-1}\sum_{j=0}^\infty 2^{-2j} \frac{\Gamma(s+2j)}{\Gamma(s-\novt+j+1)\Gamma(j+1)} \tau^{-2j},
\end{split}
\eeq
where  $d_{g_0}(w,w')$ is the distance between two points $w,w'\in \mb^{n+1}$ relative to the metric $g_0$. We should remark that the formula in  \cite{GuiZwo}  is for $\mh^n$ instead of $\mh^{n+1}$.   It is well known \cite{Hel,Laf} that $d_{g_0}$  satisfies
\beq\label{dist0}
\cosh(d_{g_0}(w,w')) =1+ \frac{2|w-w'|^2}{(1-|w|^2)(1-|w'|^2)}.
\eeq

It is proved in \cite{Gui,JosSaB},  in much greater generality, that the Schwartz kernels of $\mcp_{\bullet}(s)$ and $\mca_{\bullet}(s)$, $\bullet=0,V$,  can be obtained from $\mcr_\bullet(s,w,w')$.  We have two copies on the manifold $\mb^{n+1}$, which we denote by $\mb^{n+1}_L \times \mb^{n+1}_R$, where $R$ and $L$ stand for right and left. Let $x_R$ and $x_L$ denote the function $x$ defined in  \eqref{g0x} corresponding to each copy. For $s$ fixed, the kernel of $\mcp_{\bullet}(s)$ is a distribution $\mcp_\bullet(s,w,\theta')\in \mcd'(\mb^{n+1}\times \ms^{n})$, and similarly the kernel of the scattering matrix
$\mca_\bullet(s,\theta,\theta')\in \mcd'(\ms^n\times \ms^{n})$,
 and it turns out that  if 
$ w=(x_L,\theta)$, and $w'=(x_R,\theta')$, then for $\bullet=0,V$, the kernels of the Poisson operator $\mcp_\bullet$, its transpose $\mcp_\bullet^t$ and the scattering matrix $\mca_\bullet$ are given by
\beq\label{kernelPO}
\begin{split}
& \mcp_{\bullet}(s, w,\theta')= R_\bullet(s, w,w') x_R^{-s}\biggl|_{\{x_R=0\}}, \\
& \mcp_{\bullet}^t(s, \theta, w')= x_L^{-s} R_\bullet(s, w,w')\biggl|_{\{x_L=0\}}, \text{ and } \\
& \mca_{\bullet}(s, \theta, w')= (2s-n) x_L^{-s} R_\bullet(s, w,w') x_R^{-s}\biggl|_{\{x_L=x_R=0\}}.
\end{split}
\eeq

When $V=0$, the Schwartz kernel of the Poisson operator and the scattering matrix can be computed directly from \eqref{ker-R0}, \eqref{dist0} and \eqref{kernelPO}.  We find that, for  $x$ defined in  \eqref{g0x},
\beq\label{kernelPS}
\begin{split}
& \mcp_{0}(s,w,\theta')= c(s)  \left[ \frac{1-|w|^2}{|w-\theta'|^2}\right]^s= c(s) e^{s \bush(w,\theta)}, \;\  \mcb_0(w,\theta)= \log\left( \frac{1-|w|^2}{|w-\theta|^2}\right),\text{ and } \\
& \mca_{0}(s,\theta,\theta')= (2s-n) 4^s c(s)\left|\theta-\theta'\right|^{-2s}, \text{ where }  c(s)= \frac{1}{2\pi^\novt} \frac{\Gamma(s)}{\Gamma(s-\novt+1)}.
\end{split}
\eeq


\section{Mapping Properties of the Resolvent and the Poisson Operator}  

 In what follows we shall say that
\[
 f\in L^p(\mhn) \;\ \text{ if }  \;\ ||f||_{L^p(\mhn)}^p= \int_{\mbb} |f(w)|^p \frac{2^{n+1}dw}{(1-|w|^2)^{n+1}}<\infty, \;\ 0<p<\infty.
 \]

We will need the following estimate on the operator norm of $R_0(s):$
\begin{prop}\label{est-R0} (Proposition 3.2 of \cite{Gui-C1})  If $x$ is a boundary defining function of $\mh$, then 
\beq\label{est-R1}
 ||x^\ha R_0(\novt +i\la ) x^\ha||_{\mcl(L^2(\mhn))} \leq C |\la|^{-1}, \text{ for }  \la \in \mr \text{ and } |\la|>1.
\eeq
\end{prop}

Here and throughout the paper we shall use the notation $||T||_{\mcl(H)}$ to indicate the norm of a bounded linear operator $T:H \longmapsto H$, where $H$ is a normed vector space.  

Proposition \ref{est-R0} implies a similar estimate  for $R_V(s):$
\begin{prop}\label{est-RV}  If $V\in xL^\infty(\mbb)$, there exist $\varrho>0$ and $M>0$  such that  
\beq\label{est-resV}
 || x^\ha R_V(\novt+i\la) x^\ha||_{\mcl(L^2)} \leq  M|\la|^{-1}, \text{ provided } \la \in \mr \text{ and } |\la| > \varrho.
 \eeq
\end{prop}
\begin{proof}  Since $R_0(s)$, $s=\novt+i\la,$ is bounded in $L^2(\mhn)$, for $\im \la<0$,  the identity 
 \[
 R_0(s)(\Delta_{g_0}+V-s(n-s))= I + R_0(s) V
 \]
 is valid for $\im \la<0$, and therefore
  \[
 x^\ha R_0 (s) x^\ha \left( x^{-\ha}(\Delta_{g_0}+V-s(n-s))x^{-\ha}\right)= I + x^\ha R_0(s) x^\ha Vx^{-1}, \;\ \im \la<0.
 \]

 We know from \eqref{est-R1} and from the fact that $x^{-1} V \in L^\infty$, that 
 \[
 ||x^\ha R_0(s) x^\ha Vx^{-1}||_{\mcl(L^2(\mhn))}\leq C |\la|^{-1},  \;\ \la \in \mr, \;\  |\la|>1,
 \]
and therefore there exists $K>0$ and $\vrho>0$ such that 
\[
  ||\left(I+ x^\ha R_0(s) x^\ha Vx^{-1}\right)^{-1}||_{\mcl(L^2(\mhn))} \leq  K, \text{ provided } \la \in \mr \text{ and }  |\la|>\vrho,
  \]
 and so we have
 \[
 x^\ha R_V(s) x^\ha= \left( I + x^\ha R_0(s) x^\ha Vx^{-1}\right)^{-1} x^\ha R_0(s) x^\ha, \text{ for } \la \in \mr \text{ and }  |\la|>\vrho,
 \]
 and  therefore \eqref{est-resV} follows from \eqref{est-R1}.
  \end{proof}

\subsection{The Born Approximations of $R_V(s)$ and $\Tau_V(s)$} As in \cite{BirYaf1,BirYaf2,BulPus1} we will use the first Born approximation for $R_V(s)$ and $\Tau_V(s)$ to obtain estimates for the operator $\mce_V(s)$ defined in   \eqref{ps-rel-sm}.

One can easily check the following two identities,
 \[
 \begin{split}
 & R_V(s)= R_0(s)- R_0(s)V R_V(s), \\
& R_V(s)= R_0(s)- R_V(s)V R_0(s),
 \end{split}
 \]
 and combine them to obtain
 \beq \label{Born}
 \begin{split}
  R_V(s)- R_0(s)= - R_0(s)V R_0(s)+ R_0(s) VR_V(s)V R_0(s),
 \end{split}
 \eeq
 which is called  the first order  Born approximation for $R_V(s)$.  It follows from \eqref{Born} and \eqref{kernelPO} that
\[
\mca_V(s)-\mca_0(s)= -2i\la\bigl(\mcp_0^t(s) V \mcp_0(s)- \mcp_0^t(s) V R_V(s) V \mcp_0(s)\bigr),
\]
and so, after multiplying on the right by $\mca_0(s)^{-1}$ and using \eqref{minusL0} we obtain
\beq\label{BornSM}
\Tau_V(s)= \mca_V(s)\mca_0(s)^{-1}=\id -2i\la  \mcp_0^t(s) V \mcp_0(n-s) +2i\la \mcp_0^t(s) V R_V(s) V \mcp_0(n-s),
\eeq
which is the first order Born approximation for $\Tau_V(s)$.  We deduce that
\beq\label{mce-st}
\begin{split}
& \mce_V(s)= \Tau_V(s)-\id= \mcu_V(s) -\mcw_V(s), \text{ where } \\
& \mcu_V(s)= -2i\la  \mcp_0^t(s) V \mcp_0(n-s)  \text{ and } \\
& \mcw_V(s)=- 2i\la \mcp_0^t(s) V R_V(s) V \mcp_0(n-s).
\end{split}
\eeq
We need to obtain suitable norm estimates for the operators $\mcu_V(s)$ and $\mcw_V(s)$.

 An application of Green's identity, see for example the proof of Proposition 2.1 of \cite{Gui}, shows that
 
 \beq\label{stone0}
 R_0(s,w,w')-  R_0(n-s,w,w')= (n-2s) \int_{\ms^n} \mcp_{0}(s,w,\theta) \mcp_{0}^t(n-s, \theta,w') d\theta,
 \eeq
  On the other hand, Stone's formula  implies that the map
 \beq\label{Eins10}
 \begin{split}
 & \mcp_{0}^t(s): C_0^\infty(\mb^{n+1}) \longrightarrow C^\infty(\ms^n \times \mr_+), \\
  (\mcp_{0}^t(s) \vphi)(\theta)=&   \int_{\mb^{n+1}} \mcp_{0}^t(s,\theta,w) \vphi(w) \frac{2^{n+1}dw}{(1-|w|^2)^{n+1}}, \;\ s=\novt+i\la, \la>0,
 \end{split}
 \eeq
induces a surjective isometry
\beq\label{Eins20}
\begin{split}
\la \mcp_{0}^t(s):& \, L^2(\mh^{n+1})  \longrightarrow L^2(\mr_\la^+\times \ms^n)), \;\ s=\novt+i\la, \;\ \la\not=0, \text{ and } \\
& \int_0^\infty \int_{\ms^n} |\la \mcp_{0}^t(s) f|^2 d\theta d\la= \frac{\pi}{2} ||f||_{L^2(\mhn)},\\
\end{split}
\eeq
See for example, the proof of Proposition 2.2 of \cite{Gui}.   Moreover, it is a spectral representation of $\Delta_{g_0}$, since
\beq\label{Eins20-rep}
\mcp_{0}^t(s) \Delta_{g_0} \vphi= s(n-s) \mcp_{0}^t(s)\vphi.
\eeq

The map $\mcp_{0}^t(s)$ is the analogue of the Fourier transform on hyperbolic space, and as already mentioned, it is often called the Helgason Fourier transform.
We will  follow the strategy of \cite{BirYaf1,BulPus1} and we will need to establish weighted $L^2$ and Schatten-von Neumann  estimates  for the Poisson operator and its transpose.  We are not aware of a reference for these estimates in $\mhn$, and we prove them here.  In the Euclidean space the analogue of these estimates can be found in Chapter 8 of \cite{Yaf}.  We begin with 
\begin{prop}\label{TC-est} Let  $x$ be defined in \eqref{g0x} and let  $\mcp_{0}(s)$ be the corresponding Poisson operator defined in \eqref{Poi1}, $s=\novt+i\la$, $\la\in \mr$ and $|\la|>1$,  then the following estimates hold:
\begin{enumerate}[1.]
\item \label{prop1} For any $m\in (0,\infty)$,  there exists $C_m>0$ independent of $s$ such that
\beq\label{Eins-h1}
||\p_\la( \la \mcp_0^t(s)) x^m f||_{L^2((1,\infty)_\la \times \ms^n)} \leq C_m ||f||_{L^2(\mhn)},
\eeq 
and in particular,   there exists $C_m>0$ independent of $s$ such that 
\beq\label{est-Eins-g0}
\begin{split}
&  ||\mcp_{0}^t(s) x^m f||_{L^2(\ms^n)} \leq C_m \la^{-1} ||f||_{L^2(\mhn)}, \\
&  ||x^m \mcp_{0}  (s) \vphi ||_{L^2(\mhn)} \leq C_m \la^{-1} ||\vphi||_{L^2(\ms^n)}. 
\end{split}
\eeq
\item \label{prop2} For $m\in [1,\infty)$,  there exists $C_{m}>0$ independent of $s$  such that 
\beq\label{Eins-h11}
||\p_{\theta}( \la \mcp^t(s) x^m) f||_{L^2((1,\infty)_\la\times \ms^n)} \leq C_{m} ||f||_{L^2(\mhn)}.
\eeq 
\item \label{prop3}For $m\in[ 0,1]$ the operator $\mcp_{0}^t(s)x^m$ satisfies
\beq\label{Eins-h12}
\mcp_{0}^t(s)x^m: L^2(\mh^{n+1}) \longrightarrow H^m ((1,\infty)_\la \times \ms^n), 
\eeq
\item \label{prop4} For $m\in [\ha,\infty)$, the operator
\[
\mcp_{0}^t(s)x^m: L^2(\mh^{n+1}) \longrightarrow L^2 (\ms^n)
\]
is compact.
\end{enumerate}
\end{prop}
\begin{proof}  We know from \eqref{kernelPS} and from \eqref{metg0} that
\[
 \mcp_{0}^t(x^m f)(s,\theta)= 2^{n+1} c(s) \int_{\mb^{n+1}} e^{s \bush(w,\theta)}  x^m  f(w) \frac{dw}{(1-|w|^2)^{n+1}},
\]
 $c(s)$ and $\bush(w,\theta)$ were defined in \eqref{kernelPS}.
Therefore,
\[
\p_\la\left(\la \mcp_{0}^t(x^m f)(s,\theta)\right)= \frac{\p_\la (\la c(s))}{ \la c(s)} \la \mcp_{0}^t(x^m f)(s,\theta) + \la \mcp_{0}^t(ix^m \bush(w,\theta) f)(s,\theta).
\]

But
\[
\frac{\p_\la (\la c(s))}{ \la c(s)}= \frac{1}{\la} + \frac{\p_{\la} \Gamma(\novt+i\la)}{\Gamma(\novt+i\la)}-  \frac{\p_{\la} \Gamma(1+i\la)}{\Gamma(1+i\la)}.
\]
We recall that the quotient $\psi(z)=\frac{\Gamma'(z)}{\Gamma(z)}$ is called the digamma function, or the $\psi$ function, and it has the following asymptotic expansion
\[
\psi(z)= \log z-\frac{1}{2z} + o(\frac{1}{z}), \text{ as } z\rightarrow \infty, \;\ |\arg z|<\pi,
\]
see for example \cite{Handbook}. Hence,  
\[
\left|\frac{\p_\la (\la c(s))}{ \la c(s)}\right|\leq \frac{C}{\la},  \;\ s=\novt+i\la, \;\ \la\in \mr_+.
\]

Since
\[
|\bush(w,\theta)| \leq C| \log (1-|w|)|, \text{ it follows that }   |x^m \bush(w,\theta)| \leq C \text{ for } m> 0,
\]
and so  it follows from \eqref{Eins20} that, for $m>0$,
\[
\begin{split}
& || \la \mcp_{0}^t(x^m f)||_{L^2((1,\infty)_\la \times \ms^n)} \leq C ||f||_{L^2(\mh^{n+1})} \text{ and } \\
& ||\p_\la\left( \la \mcp_{0}^t(x^m f)\right)||_{L^2((1,\infty)_\la \times \ms^n)} \leq C ||f||_{L^2(\mh^{n+1})}.
\end{split}
\]
This  proves \eqref{Eins-h1} and it implies that $\la \mcp_{0}^t(x^m f)(\la,\theta) \in H^1((1,\infty)_\la; L^2(\ms^n))$ and in particular this implies the first estimate in \eqref{est-Eins-g0}. The second estimate follows from the first and this proves the first item.

To prove the second statement, pick coordinates $\theta= (\theta',\theta_n)$, $\theta_n^2=1-|\theta'|^2$, $|\theta_n|>\eps$.
\[
\p_{\theta_j}\mcp_{0}^t(x^m f)(s,\theta)=s \mcp_{0}( x^m \p_{\theta_j} \mcb(w,\theta) f), \;\ j<n,
\]
but 
\[
\p_{\theta_j} \bush(w,\theta)= 2\frac{w_j-\theta_j}{|w-\theta|^2}+ 2\frac{w_n-\theta_n}{|w-\theta|^2}\frac{\theta_j}{\theta_n},
\]
and since $|\theta_n|>\eps$,  $|\p_{\theta_j} \bush(w,\theta)| \leq \frac{2}{1-|w|}$, it follows that $|x^m \p_{\theta_j} \bush(w,\theta)|\leq C x^{m-1}$, provided $m>1$. Since $\ms^n$ can be covered by finitely many such charts,  \eqref{Eins-h11} follows from \eqref{Eins20}, and this proves \eqref{Eins-h11}.

We have shown that if $m\geq 1$, $\la \mcp_0(s) x^m f \in H^1((1,\infty)_\la \times  \ms^n)$.   Now notice that for $x>0$ and $a>0$,
\[
\begin{split}
& x^{a+ib}= x^a( \cos(b \log x)+ i\sin(b\log x)), \\ 
& \p_a x^{a+ib}= x^{a+ib} \log x \text{ and } \p_b x^{a+ib}= ix^{a+ib} \log x,
\end{split}
\]
and so these are continuous functions up to $x=0$ and so if $\re z>0$,  the function
\[
f(x,z)= \left\{\begin{array}{ll} x^z  \text{ if } x>0, \\ 0 \text{ if } x=0 \end{array}\right.
\]
is holomorphic in $z$ (for fixed $x$)  and continuous in $x$ (for fixed $z$),  provided $\re z>0$.

We consider the family of operators
\[
\mcf(\alpha)= s\mcp_0^t(s) x^\alpha, \;\  \re \alpha > 0.
\]
The family $\mcf(\alpha)$ is holomorphic in $\alpha$ for $\re\alpha>0$.  Also notice that 
 it follows from  \eqref{Eins20} that
\beq\label{aux1}
||s\mcp_0^t(s) x^{ib} f||_{L^2((1,\infty)_\la \times \msn)} \leq C ||x^{ib} f||_{L^2(\mhn)}= C||f||_{L^2(\mhn)}.
\eeq
This implies that the family $\mcf(\alpha)$ is holomorphic for $\re\alpha>0$ and continuous up to $\re\alpha=0$. We also know from
 \eqref{Eins-h1}, \eqref{Eins-h11} and \eqref{aux1}  that
\[
\begin{split}
& ||\mcf(\alpha) f||_{L^2( (1,\infty)_\la \times \msn)} \leq C ||f||_{L^2(\mhn)}, \text{ provided } \;\ \re\alpha=0, \\
& ||\mcf(\alpha) f||_{H^1( (1,\infty)_\la \times \msn)} \leq C ||f||_{L^2(\mhn)}, \text{ provided } \;\ \re\alpha\geq 1.
\end{split}
\]
It follows from Stein's interpolation theorem \cite{Stein}, that for $m\in (0,1)$,
\[
||\mcf(\alpha) f||_{H^m( (1,\infty)_\la \times \msn)} \leq C ||f||_{L^2(\mhn)}, \text{ provided } \;\ \re\alpha=m,
\]
and this proves \eqref{Eins-h12}. 

Since $\ms^n$ is compact,  if $s=\novt+i\la$, $|\la|>1$,  is fixed and $m>\ha$, it follows from \eqref{Eins-h12} and Sobolev restriction theorem that
\[
||\mcp^t_0(s) x^m f||_{H^{m-\ha}(\msn)} \leq C ||f||_{L^2(\mhn)},
\]
 and  since $m-\ha>0$, we conclude  from the Sobolev embedding theorem that  $\mcp_0^t(s) x^m$, as an operator acting on $L^2(\msn)$,  is compact for $m>\ha$.   
 To show that the result still holds for $m=\re\alpha=\ha$, notice that
 \[
 \mcp_0^t(s) x^m f-\mcp_0^t(s) x^\ha f=  2^{n+1}c(s) \int_{\mb^{n+1}} e^{s \bush(w,\theta)}  (x^m- x^\ha)  f(w) \frac{dw}{(1-|w|^2)^{n+1}}.
 \]
But we notice that
 \[
 \begin{split}
 & x^m-x^\ha = x^{\ha} \bigl(x^{m-\ha}-1\bigr)=(m-\ha)x^{\ha}(\log x) q(m,x), \;\  q(m,x)\defn\frac{x^{m-\ha}-1}{(m-\ha) \log x},
 \end{split}
 \]
 and one can check that $0<q(m,x)<1$, for $x\in (0,1)$ and $m>\ha$.  It follows that, if $m\in (\ha,1)$,
 \[
  \mcp_0^t(s) x^m f-\mcp_0^t(s) x^\ha f=  2^{n+1}c(s) \int_{\mb^{n+1}} e^{s \bush(w,\theta)} x^{\ha} (\log x) (m-\ha) q(m,x)  f(w) \frac{dw}{(1-|w|^2)^{n+1}}, 
  \]
 and it follows from \eqref{est-Eins-g0} that
 \[
 || \mcp_0^t(s) x^m f-\mcp_0^t(s) x^\ha f||_{L^2(\msn)} \leq C \biggl(m-\ha\biggr)\;  ||f||_{L^2(\mhn)}, \;\ m\in (\ha,1).
 \]
 
 This means that the sequence of compact operators $\mcp_0^t(s) x^m$, $m>\ha$,  converges strongly to $\mcp_0^t(s) x^\ha$ as $m\rightarrow \ha$ and so $\mcp_0^t(s) x^\ha$  is compact , and this ends the proof of the  Proposition.
\end{proof}

Next we obtain the following estimates for $\mce_V(s):$
\begin{prop}\label{est-EV}  If  $m\in (1,\infty)$ and $V\in x^mL^\infty(\mhn)$ is real valued,  then there exist $R>0$ and $M>0$ such that
\beq\label{est-EV1}
||\Tau_V\bigl(\novt+i\la\bigr)-\id||_{\mcl(L^2(\msn))}\leq M||x^{-m}V||_{{}_{L^\infty}}|\la|^{-1}, \text{ provided } \la \in \mr, \;\ |\la| >R.
\eeq
\end{prop}
\begin{proof}  If $s=\novt+i\la,$ $\la \in \mr,$ we deduce from \eqref{est-Eins-g0} that there exists $M>0$, independent of $\la$  such that if $|\la|>1,$
\beq\label{est-P1}
||\la \mcp_0^t(s) V \mcp_0(n-s)||_{\mcl(L^2(\msn)}=||\la \mcp_0^t(s)x^{\frac{m}{2}} (x^{-m} V) x^{\frac{m}{2}} \mcp_0(n-s)||_{\mcl(L^2(\msn)}
\leq M ||x^{-m}V||_{{}_{L^\infty}}|\la|^{-1}, 
\eeq
  Similarly, we combine \eqref{est-Eins-g0} and \eqref{est-resV} to deduce that there exist $R>0$ and $M>0$, independent of $\la$, such that
\beq\label{est-P2}
|| \la \mcp_0^t(s) V R_V(s) V \mcp_0(n-s)||_{\mcl(L^2(\msn))} \leq M||x^{-m}V||_{{}_{L^\infty}} |\la|^{-2} \text{ for } |\la| >R.
\eeq
\end{proof}

\subsection{Schatten-von Neumann Estimates for $\Tau_V(s)$}  
We recall some definitions and results about this topi and we refer the reader to \cite{GolKre,Sim} for more details. If  $H$ and $\mch$ are Hilbert spaces, and $T: H\longrightarrow \mch$ is a linear bounded compact operator, then $T^*T: H \longrightarrow H$ is a compact self-adjoint nonnegative operator, and one defines $|T|=\sqrt{T^*T}$. This is also compact, nonnegative and self-adjoint.  Let $s_j(T)$ denote the eigenvalues of the operator $|T|, $ which are called the singular values of $T$. They are nonnegative and arranged in such a way that $s_j(T)\geq s_{j+1}(T)$, and counted with multiplicity, and  converge to zero as $j$ goes to infinity.

 The Schatten-von Neumann classes $\Sigma_p(H,\mch)$,  $p\geq 1$, are defined to be the space of compact operators 

\beq\label{scat1}
\begin{split}
& T: H \longrightarrow \mch  \text{ such that } \\
 ||T||_{{}_{\Sigma_p}}^p =\sum_{j=1}^\infty  & s_j(T)^p<\infty, \text{ provided } 1 \leq p < \infty. \\
 \text{ If } p=\infty,  \text{ in which case }  &T \text{ is just compact, }  ||T||_{{}_{\Sigma_\infty}}=||T||, \text{ the norm of } T.
\end{split}
\eeq

This defines a norm, and these spaces are nested 
\[
 \Sigma_p(H,\mch)\subset \Sigma_q(H,\mch) \text{ and }  ||T||_{{}_{\Sigma_q}}\leq  ||T||_{{}_{\Sigma_p}}, \text{ provided } p \leq q,
 \]
 and  form ideals in the space of bounded operators,
\[
\begin{split}
 A \in\  & \Sigma_p(H,\mch_1)  \text{ and }   B \in \mcl(\mch_1,\mch_2) \text{ is a bounded operator with norm } ||B||, \\
 & \text{ then } BA \in \Sigma_p(H,\mch_2) \text{ and } 
   ||BA||_{{}_{\Sigma_p}}\leq   ||A||_{{}_{\Sigma_p}} ||B||.
\end{split}
\]

 These spaces also satisfy H\"older's inequality:
\beq\label{Holder}
\begin{split}
& \text{ if } A \in  \Sigma_{p_1}(H,\mch_1) \text{ and } B \in \Sigma_{p_2}(\mch_1, \mch_2), \;\ 1< p_j\leq \infty,  \text{ then } \\
& BA \in \Sigma_p(H,\mch_2) \text{ and } ||BA||_{\Sigma_p} \leq ||A||_{\Sigma_{p_1}} ||B||_{\Sigma_{p_2}}, \\
& \text{ provided } \frac{1}{p}=\frac{1}{p_1}+ \frac{1}{p_2}.
\end{split}
\eeq

$\Sigma_2$ is also called the class of Hilbert-Schmidt operators and if $T:H \longrightarrow \mch$, and $\{\vphi_j, j \in \mn\}$ is an orthonormal basis for $H$,
\[
||T||_{{}_{\Sigma_2}}=\sum_j ||T \vphi_j||_{\mch}^2,  \text{ and is independent of the choice of the basis.}
\]
If $T: L^2(X) \longrightarrow L^2(Y)$ is defined by 
\[
Tf (y)= \int_{X} K_T(y,z) f(z) \ dz,
\]
then
\[
||T||_{{}_{\Sigma_2}}= \int_{Y\times X} |K_T(y,z)|^2 dz dy.
\]

  If $\mch=H$, the class $\Sigma_1(H)$ consists of trace class operators, and the trace of the operator $T$ is given by
\[
Tr(T)= \sum_{n=1}^\infty \lan T f_n,f_n\ran_{H},  \text{ where } \lan, \ran \text{ is the inner product of } H \text{ and } f_n \text{ is a basis for } H,
\]
which converges absolutely and is independently of the choice of the basis.

We will need the following Schatten-von Neumann norm estimates for  the weighted Poisson operator:
\begin{prop}\label{Trace-est} Let $x$ be defined in \eqref{g0x} and let $\mcp_{0}(s)$ be defined in \eqref{Poi1},  then the following are true for $s=\novt+i\la$, provided $\la \in \mr$ and $|\la|>1:$ 
\begin{enumerate}[1.]
\item \label{Tr1} Let $\mck_m(s,\theta,w)\in \mcd'(\msn\times \mb^{n+1})$ denote the Schwartz kernel of the operator $\mcp_{0}^t(s) x^m$. Then, for $m\in (n, \infty)$,  there exists a constant $C_m$ independent of $s$ such that
\beq\label{HS-g0}
\begin{split}
 ||\mck_m(s,\theta,w)||_{L^2( \msn\times \mhn)} \leq   C_m |\la|^{\frac{n}{2}-1},
\end{split}
\eeq
and hence $\mcp_{0}^t(s) x^m \in \Sigma_2(L^2(\mhn), L^2(\msn))$ and its Hilbert-Schmidt norm  satisfies
\beq\label{HS1-g0}
||\mcp_{0}^t(s) x^m||_{{}_{\Sigma_2(L^2(\mhn),L^2(\msn))}} \leq C_m|\la|^{\novt-1}.
\eeq
\item \label{Tr2} If  $m\in (\ha,\infty)$, $p\in [2, \infty)$, then $\mcp^t_0(s) x^m \in \Sigma_p(L^2(\mhn),L^2(\msn))$  and 
\beq\label{three-line}
||\mcp_0^t(s) x^m||{{}_{{}_{\Sigma_p(L^2(\mhn), L^2(\msn))}}}\leq C |\la|^{\frac{n}{p}-1}, \text{ provided }  p> \frac{2(n-\ha)}{m-\ha}.
\eeq
\item \label{Tr3} If $m\in (1,\infty)$, $p\in[ 1,\infty)$ and $p> \frac{2(n-\ha)}{m-1}$,  then $ \mcp_0^t(s) x^m \mcp_0(n-s)\in \Sigma_p(L^2(\msn),L^2(\msn))$ and 
\beq\label{three-line1}
\begin{split}
& ||\la \mcp_0^t(s) x^m \mcp_0(n-s)||_{{}_{\Sigma_p(L^2(\msn), L^2(\msn))}}\leq C |\la|^{\frac{n}{p}-1}, 
\end{split}
\eeq
\item If $m\in (1,\infty)$ and  $p\in[ 1,\infty)$ and $p> \frac{2(n-\ha)}{m-1}$, then there exists $\rho>0$ such that for $ |\la|>\rho$ we have
\beq\label{three-line2}
\begin{split}
&  \mcp_0^t(s) x^m R_V(s) x^m\mcp_0(n-s)\in \Sigma_p(L^2(\msn),L^2(\msn)), \text{ and } \\ 
& ||\la \mcp_0^t(s) x^m R_V(s) x^m \mcp_0(n-s)||_{{}_{\Sigma_p(L^2(\msn), L^2(\msn))}}\leq C |\la|^{\frac{n}{p}-2},
\end{split}
\eeq
\end{enumerate}
\end{prop}
\begin{proof} We know from \eqref{kernelPS} that
\[
\mck_m(s,w,\theta)= c(s) \left(\frac{1-|w|^2}{|w-\theta|^2}\right)^s x^m.
\]
Since  $|w-\theta|^2\geq (1-|w|)^2$,  $s=\novt+i\la$ and $x=\frac{1-|w|}{1+|w|}$, we have 
\[
\begin{split}
& ||\mck_m(s,z,\theta) ||_{L^2(\mhn\times \msn)}^2 \leq \\
&2^{n+1} |c(s)|^2 \int_{\ms^{n}} \int_{\mb^{n+1}}  \frac{(1-|w|^2)^n}{|w-\theta|^{2n}} \frac{(1-|w|)^{2m}}{(1+|w|)^{2m}} \frac{dw}{(1-|w|^2)^{n+1}} d\theta  dw\leq \\
& M |c(s)|^2 \int_{\mb^{n+1}} (1-|w|)^{2m-2n-1} dw \leq M' |c(s)|^2, \text{ provided } m> n.
\end{split}
\]
We know from \cite{ErdTri} that
\beq\label{Tricomi}
\frac{\Gamma(z+\alpha)}{\Gamma(z+\beta)}= z^{\alpha-\beta}\left(1+ \frac{(\alpha-\beta)(\alpha+\beta-1)}{2z} + O(|z|^{-2})\right) \text{ as } |z|\rightarrow \infty,
\eeq
and in particular,
\beq\label{quotient}
c(\novt+i\la)= \frac{1}{2\pi^n} \la^{\novt-1}( 1+ O(\frac{1}{|\la|})), \text{ as } \la\rightarrow \infty,
\eeq
and so for $s=\novt+i\la$,
\[
|c(s)|^2 \leq C \la^{n-2}, \;\ \la >>1.
\]
This proves \eqref{HS-g0}, and \eqref{HS1-g0} follows from it.

To prove the second item, we proceed as in  as in \cite{Yaf}.    For fixed $s$, we consider the family
\[
\mcf(\alpha) = \mcp_{0}^t(s) x^\alpha, \;\ \alpha \in \CC, \;\ \re\alpha>0,
\]
which is holomorphic in $\alpha$. We have shown in Proposition \ref{TC-est} and in \eqref{HS-g0} that
\begin{enumerate}[1.]
\item $\mcf(\alpha)\in \Sigma_\infty(L^2(\mhn), L^2(\msn))$ and $||\mcf(\alpha)||\leq C\la^{-1}$ if $\re\alpha\geq \ha$
\item $\mcf(\alpha) \in \Sigma_2(L^2(\mhn), L^2(\msn))$  and $||\mcf(\alpha)||_{{}_{\Sigma_2}} \leq C \la^{\novt-1}$ if $\re\alpha> n$, 
\end{enumerate}
and so it follows from the three-line theorem for operator valued functions in Schatten classes, see Theorem 2.6 from  Chapter 0 of \cite{Yaf} that for all  $m_2>n$ and $m \in(\ha,m_2)$, and $\re \alpha=m$,
\[
\begin{split}
\mcp_{0}^t(s) x^\alpha  \in \Sigma_p(& L^2(\mhn),  L^2(\msn)), \text{ provided } p= \frac{2(m_2-\ha)}{m-\ha}, \text{ and }\\
&  \text{ and }  ||\mcf(\alpha)|_{{}_{\Sigma_p}}\leq C \la^{\frac{n}{p}-1}.
\end{split}
\]
Since $m_2>n$, this proves \eqref{three-line}.

 To prove \eqref{three-line1} we use that \eqref{three-line} holds for $m>1$, and for $p_1,p_2\geq 2$ such that  $p_1,p_2>\frac{2(n-\ha)}{\frac{m}{2}-\ha}$. Therefore, if $p$ is such that $\frac{1}{p}=\frac{1}{p_1}+\frac{1}{p_2}$, then  $p\geq 1$ and $ \frac{1}{p}< \frac{m-1}{2(n-\ha)}$.  So we  have from \eqref{Holder}  and \eqref{three-line} that 
 \[
 ||i\la \mcp_0^t(s) x^m \mcp_0(n-s)||{{}_{\Sigma_p}}\leq  ||\la\mcp_0^t(s) x^{\frac{m}{2}}||_{{}_{\Sigma_{p_1}}} ||x^{\frac{m}{2}}\mcp_0(n-s)||_{{}_{\Sigma_{p_2}}} \leq C 
 |\la|^{\frac{n}{p_1}} |\la|^{\frac{n}{p_2}-1}  \leq C |\la|^{\frac{n}{p}-1} 
 \]
 Similarly, to prove \eqref{three-line2} we write
 \[
  i\la \mcp_0^t(s) x^m R_V(s) x^m\mcp_0(n-s)= i\la \mcp_0^t(s) x^{m-\ha} x^\ha R_V(s) x^\ha x^{m-\ha}\mcp_0(n-s).
  \]
  Again, we use \eqref{Holder} and \eqref{three-line} and \eqref{est-resV} to show that, provided $p\geq 1$, $p>\frac{n-\ha}{m-1}$, we have
  \[
  \begin{split}
  & || i\la \mcp_0^t(s) x^{m-\ha} x^\ha R_V(s) x^\ha x^{m-\ha}\mcp_0(n-s)||_{{}_{\Sigma_{p}}}\leq \\
   &  || i\la \mcp_0^t(s) x^{m-\ha} ||_{{}_{\Sigma_{p_1}}} \ || x^\ha R_V(s) x^\ha|| \ ||x^{m-\ha}\mcp_0(n-s)||_{{}_{\Sigma_{p_2}}} \leq C |\la|^{\frac{n}{p_1}+\frac{n}{p_2}-2}=
    C |\la|^{\frac{n}{p}-2}.
    \end{split}
    \]
    
  Therefore, \eqref{three-line2} follows from \eqref{three-line}
  and \eqref{est-Eins-g0}. This ends the proof of the  Proposition.
\end{proof}

We apply these results to obtain the following Schatten-von Neumann estimates for $\Tau_V(s)-\id$, $\mcu_V(s)$ and $\mcw_V(s)$ defined in \eqref{mce-st}:

\begin{prop}\label{est-mcu} If $V \in x^mL^\infty(\mb^{n+1})$, $m\in (1,\infty)$, there exists $\rho>0$ such that the  operators
$\mcu_V(s)$ and $\mcw_V(s)$, $s=\novt+i\la,$ $\la\in \mr$, defined in \eqref{mce-st} satisfy $\mcu_V(s) , \mcw_V(s) \in \Sigma_p(L^2(\msn), L^2(\msn))$,  for any $p\geq 1$ such that $p>\frac{2(n-\ha)}{m-1}$ and  $|\la|>\rho$. Moreover, 
\beq\label{mcu}
\begin{split}
& ||\mcu_V(s)||_{{}_{\Sigma_p}} \leq C ||x^{-m}V||_{{}_{L^\infty}} |\la|^{\frac{n}{p}-1} \text{ and } 
   ||\mcw_V(s)||_{{}_{\Sigma_p}} \leq C ||x^{-m}V||_{{}_{L^\infty}} |\la|^{\frac{n}{p}-2}.
\end{split}
\eeq
These can also be stated in terms of  $\Tau_V(s)-\id\in \Sigma_p(L^2(\msn), L^2(\msn))$ and 
\beq\label{Tau}
\begin{split}
& ||\Tau_V(s)-\id|||_{{}_{\Sigma_p}} \leq C ||x^{-m}V||_{{}_{L^\infty}}|\la|^{\frac{n}{p}-1} \text{ and } \\
& ||\Tau_V(s)-\id- \mcu_V(s)|||_{{}_{\Sigma_p}} \leq C  ||x^{-m}V||_{{}_{L^\infty}}|\la|^{\frac{n}{p}-2},
\end{split}
\eeq
 provided $p\geq 1$, $p>\frac{2(n-\ha)}{m-1}$ and  $|\la|>\rho$.
\end{prop}
\begin{proof} The proof follows from the simple observation that  for $W=x^{-m}V\in L^\infty(\mhn)$, $\mcu_V(s)$ satisfies
\[
\mcu_V(s)= 2i\la \mcp_0(s)^tV \mcp_0(n-s)=  2i\la \mcp_0(s)^t x^{\frac{m}{2}} W(z) x^{\frac{m}{2}} \mcp_0(n-s), \;\ W\in L^\infty,
\]
and since $m>1$, one can pick $p_1\geq 2$ and $p_2\geq 2$ such that $p_j>\frac{2(n-\ha)}{m-\ha}$ and $\frac{1}{p}=\frac{1}{p_1}+\frac{1}{p_2}$ and in view of \eqref{three-line}
\[
||\mcu_V(s)||_{\Sigma_p} \leq 2 ||\la \mcp_0(s)^t x^{\frac{m}{2}}||_{\Sigma_{p_1}} ||W||_{{}_{L^\infty}} ||x^{\frac{m}{2}} \mcp_0(n-s)||_{\Sigma_{p_2}} \leq C |\la|^{\frac{n}{p}-1}.
\]
 Similarly,
\[
\begin{split}
& \mcw_V(s)= 2i\la \mcp_0(s)^tV R_V(s) V \mcp_0(n-s)=  2i\la \mcp_0(s)^t x^{m-\ha} W(z) x^\ha R_V(s) x^{\ha} W(z)  x^{m-\ha}\mcp_0(n-s), 
\end{split}
\]
and since $m-\ha>\ha$ the same argument applies, fo $|\la|>\rho$ such that \eqref{est-resV} holds.
\end{proof}

 As a consequence of this we obtain the following, which essentially is Lemma 3.3 of \cite{BulPus1}:
\begin{prop}\label{eig-ap} Let $m\in (1,\infty)$  and $V\in x^m L^\infty(\mhn)$,  be real valued. There exist $R>0$ and $C>0$ such that if $p\in \mn$, $p\geq 1$, and  $p>\frac{2(n-\ha)}{m-1}$,  $s=\novt+i\la$, $\la \in \mr$ and $|\la|>R$, then 
\beq\label{eig-ap0}
\biggl|  \sum_{e^{i\del_j(s)}\in \spec \Tau_V} \bigl(\la\del_j(s)\bigr)^p-Tr\bigl(\la \im \Tau_V(s)\bigr)^p \biggr|\leq C ||x^{-m}V||_{L^\infty}|\la|^{n-2}.
\eeq
\end{prop}  
\begin{proof} 
Since $V$ is real valued, $\im \Tau_V(s)=\im(\Tau_V-\id)$ is compact and self-adjoint,
\[
Tr\bigl(\la \im \Tau_V(s)\bigr)^p=\sum_{e^{i\del_j(s)}\in \spec \Tau_V} \bigl(\la \sin(\del_j(s))\bigr)^p,
\]
 and so \eqref{eig-ap0} is equivalent to
 \[
\mathfrak{S}(s)\defn \biggl| \sum_{e^{i\del_j(s)}\in \spec \Tau_V} \bigl(\la\del_j(s)\bigr)^p- \sum_{e^{i\del_j(s)}\in \spec \Tau_V} \bigl(\la \sin(\del_j(s))\bigr)^p \biggr| \leq C ||x^{-m}V||_{L^\infty} |\la|^{n-2}.
\]
In view of \eqref{est-EV1}, there exists $R_0>0$ such that  $|\sin \del_j(s)|< C \| x^{-m}V\|_{L^\infty} |\la|^{-1}$,  for all $j$, provided $|\la|>R_0.$ In particular this implies that there exists $C>0$ such that, for $R$ large enough, 

\[
\ha|\del_j(s)| \leq |\sin \del_j(s)| \leq |\del_j(s)| \text{ and }  |\del_j(s)-\sin(\del_j(s))|  \leq  C|\sin(\del_j(s))|^3, \;\ |\la|>R.
\]
Therefore
\[
\begin{split}
& \bigl|(\la \del_j(s))^p-(\la \sin(\del_j(s))^p\bigr|= |\la|^p\bigl| (\del_j(s))^p- (\sin \del_j(s))^p\bigr|=\\
&|\la|^p |\del_j(s)- \sin(\del_j(s))| Q(|\del_j(s)|, |\sin(\del_j(s))|),
\end{split}
\]
where $Q(x,y)$ is a homogeneous polynomial of degree $p-1$, and therefore,
\[
\bigl|\la \del_j(s))^p-(\la \sin(\del_j(s))^p\bigr|\leq C |\la|^p|\bigl|\sin( \del_j(s))^{p+2}\bigr|.
\]
Again using that  $\im\Tau_V(s)$ is compact and self-adjoint, 
\[
\sum_{j} \bigl| \sin( \del_j(s))^{p+2}\bigr|= ||\im \Tau_V||_{\Sigma_{p+2}}^{p+2}.
\]
So we conclude that, in view of \eqref{Tau}, 
\[
\mathfrak{S}(h) \leq C |\la|^p  ||\im \Tau_V||_{\Sigma_{p+2}}^{p+2} \leq C ||x^{-m}V||_{L^\infty} |\la|^{n-2},
\]
and this ends the proof of the proposition.
\end{proof}

We remark that it follows from \eqref{kernelPS}, and the fact that $V$ is real valued, that 
\beq\label{imagin}
\mcu_V^*(s)= -\mcu_V(s) \text{ and so } \im \mcu_V(s) = \frac{1}{2i} (\mcu_V(s)-\mcu_V(s)^*)= \frac{1}{i} \mcu_V(s),
\eeq
and as in Lemma 2.2 of \cite{BulPus1}, we consider the operator
\[
\im \Tau_V(s)= \frac{1}{2i}(\Tau_V(s)-\Tau_V^*(s))= \frac{1}{i} \mcu_V(s)+ \im \mcw_V(s),
\]
and we  have 
\begin{prop}\label{ptr} Let $m\in (1,\infty)$ and let  $V \in x^m L^\infty(\mbb)$ be real valued.  If  $\rho$ is as in Proposition \ref{est-mcu},  then for any $p\in \mn$, such that $p\geq 1$ and $p>\frac{2(n-\ha)}{m-1}$ and $s=\novt+i\la$, with $\la \in \mr$ and $|\la|>\rho$, then 
\beq\label{ptr1}
| \tr\left[ \left(\la \im\Tau_V(s)\right)^p - (\la \mcu_V(s))^p\right]| \leq C  ||x^{-m}V||_{{}_{L^\infty}}|\la|^{n-1}.
\eeq
\end{prop}
\begin{proof} We present the proof for the convenience of the reader, but we just follow the arguments of \cite{BulPus1}.  First we realize that if $A=\la\im \Tau_V(s)$ and $B=\frac{\la}{i}\mcu_V(s)$ then
\beq\label{amb1}
A^p-B^p= \sum_{j=0}^{p-1} A^j(A-B) B^{p-1-j},
\eeq
Therefore, by \eqref{Holder}
\beq\label{amb2}
||A^p-B^p||_{{}_{\Sigma_1}} \leq \sum_{j=0}^{p-1} ||A||_{{}_{\Sigma_p}}^j \ ||A-B||_{{}_{\Sigma_p}} \  ||B||_{{}_{\Sigma_p}}^{p-j-1}.
\eeq
But we know from Proposition \ref{est-mcu} that
\[
\begin{split}
& ||A||_{{}_{\Sigma_p}}  \leq C  ||x^{-m} V||_{{}_{L^\infty}} |\la|^{\frac{n}{p}}, \\
& ||B||_{{}_{\Sigma_p}}  \leq C ||x^{-m}V ||_{{}_{L^\infty}} |\la|^{\frac{n}{p}} \text{ and } \\
& ||A-B||_{{}_{\Sigma_p}}  \leq C  ||x^{-m} V||_{{}_{L^\infty}}|\la|^{\frac{n}{p}-1},
\end{split}
\]
and this gives \eqref{ptr1}.
\end{proof}

As a consequence of Propositions \ref{ptr} and \ref{eig-ap} we obtain
\begin{prop}\label{corol} Let $m\in (1,\infty)$ and let  $V\in x^m L^\infty(\mbb)$ be real valued. There exists $R>0$ such that if  $p\in \mn$, $p\geq 1$ and $p>\frac{2(n-\ha)}{m-1}$, and $s=\novt+i\la$, $\la \in \mr$ and $|\la|>R$, then 
\beq\label{ptr1N}
\biggl| \sum_{e^{i\del(s)}\in \spec \Tau_V} \bigl(\la\del_j(s)\bigr)^p - Tr (\la \mcu_V(s))^p\biggr| \leq C  ||x^{-m}V||_{{}_{L^\infty}}|\la|^{n-1}.
\eeq

\end{prop}

\section{The structure of the operator $\mcu_V$ for $V\in C_0^\infty(\mb^{n+1})$}  \label{SCP}

We will show that if $V\in C_0^\infty(\mb^{n+1})$, then $\mcu_V(\novt+\ioh)$ is in a very special class of semiclassical pseudodifferential operators and in Section \ref{TRP} we will use the calculus of such operators to prove the trace formulas.   We should remark that this result is  a particular case of a more general theorem. Notice that the  Poisson operator for $\mhn$ is given by \eqref{kernelPS} and so
\[
( V\mcp_0 f)(h,w)=  c\bigl(\novt+\ioh\bigr) \int_{\msn}  V(w) e^{ (\novt+\ioh) \mcb_0(w,\theta)} f(\theta) \ d\theta, \;\ 
\mcb_0(w,\theta)= \log\biggl( \frac{1 -|w|^2}{|w-\theta|^2 }\biggr).
\]

If $V\in C_0^\infty,$ then  $|w|<1-\eps$ for $\eps>0$, and the phase $\mcb_0(w,\theta)$  is a $C^\infty$ function and so this is a semiclassical Fourier integral operator, see for example \cite{Dui1,GuiSte}, associated with the Lagrangian submanifold $\La$ of $T^*(\mr^{n+1} \times  \msn)$ defined by
\[
\La= \bigl\{(w,\varsigma, \theta, \eta): \varsigma = \p_w \mcb_0(w,\theta) \text{ and }  \eta = \p_\theta \mcb_0(w,\theta)\bigr\}.
\]
  So $\mcu_V(\novt+\ioh)$ is the composition of two Fourier integral operators, and the result proved in this section can be thought of as an application of  the  composition theorem for semiclassical Fourier integral operators, see for example \cite{GuiSte}, or \cite{Hor} for regular Fourier integral operators.  In this particular case we have an explicit formula and we can prove this result directly using stationary phase methods (which is how the general case is proved),  and we start by obtaining an asymptotic expansion in $h$ for the  Schwartz kernel of $\ooh \mcu_V(\novt+\ioh)$.  

\begin{lemma}\label{princ1}  If $V\in C_0^\infty(\mhn)$ and $\del>0$, the Schwartz kernel of the operator $\ooh \mcu_V(\novt+\ioh)$  defined in \eqref{mce-st}  satisfies
\beq\label{princ2}
\begin{split}
  & \frac{1}{h}  \mcu_V(\novt  +\frac{i}{h}, \theta,\theta')=  O(h^\infty) \text{ if }  |\theta-\theta'|> \del,  \text{ and otherwise we have }  \\
  & \frac{1}{h}  \mcu_V(\novt  +\frac{i}{h},  \theta,\theta') =\\
  &   \frac{2^n}{ (\pi h)^n}\int_{\mr^{n+1}}   e^{-\frac{2i}{h}\lan y, \theta-\theta' \ran} \bigl(H(\frac{y}{|y|^2},\theta,\theta') + hR(h,y,\theta,\theta')\bigr) |y|^{-2(n+1)} \ dy + O(h^\infty),  \text{where }  \\
  & H(z,\theta,\theta')=   (|z-\ga|\,  |z+\ga|)^{-n} (1-|z+\ha(\theta+\theta')|^2)^{-1}V(z+\ha(\theta+\theta')), \;\ \ga=\ha(\theta-\theta')  \text{ and } \\
   &  R(h,y,\theta,\theta' )   \sim  \sum_{j=0}^\infty h^j  \ R_j(y,\theta,\theta'), \;\ R_j\in C_0^\infty(\mr^{n+1}\times \msn\times \msn),  \text{ with support independent of } j
\end{split}
\eeq
Here the notation $O(h^\infty)$ means that it  is an integral of the form
\[
\begin{split}
& O(h^\infty)= \int_{\mr^{n+1}}   e^{-\frac{2i}{h}\lan y, \theta-\theta' \ran}  a(h,y,\theta,\theta') \ dy \text{ for some } a(h,y,\theta,\theta') \\
&  \text{  which is compactly supported and }  |(\p_y, \p_\theta,\p_\theta')^\alpha a| \leq C_{\alpha,N} h^{N} \text{ for all } N, \alpha.
\end{split}
\]
\end{lemma}

\begin{proof}   It follows from \eqref{kernelPS}  and \eqref{mce-st} that the kernel of $\mcu_V(s)$, $s=\novt+i\la$,  $\la \in \mr\setminus 0,$ satisfies
\[
\begin{split}
 \la \mcu_V(s)(\theta,\theta')= &  2^{n+2}\la^2 c(s)c(n-s) \int_{\mbb} e^{i\la \left(\log |w-\theta'|^2- \log |w-\theta|^2\right)} G(w,\theta,\theta')dw, \text{ where } \\
& G(w,\theta,\theta')=  (|w-\theta|\, |w-\theta'|)^{-n} (1-|w|^2)^{-1}V(w).
\end{split}
\]
Since $V\in C_0^\infty(\mbb)$, and $\theta\in \msn$, there exists $\eps>0$ such that $|w-\theta|>\eps>0$,  for all $\theta\in \msn$, provided $w \in \supp(V)$. We set $w=z+\ha(\theta+\theta')$ and so 
\beq\label{def-Z}
\begin{split}
& w-\theta= z-\ga \text{ and } w-\theta'= z+\ga, \;\ \ga=\frac{\theta-\theta'}{2}, \\
&  |z-\ga|>\eps >0 \text{ and } |z+\ga|>\eps >0.
\end{split}
\eeq

As usual we write $h=\la^{-1}$, and so we have
\beq\label{defH}
\begin{split}
& \la \mcu_V(s)(\theta,\theta')= \frac{2^{n+2}}{h^2} c(\novt+\frac{i}{h})c(\novt-\frac{i}{h}) \Upsilon_V(h,\theta,\theta'), \\
& \Upsilon_V(h,\theta,\theta')=\int_{\mr^{n+1}} e^{\ioh  \vphi(z,\ga)} H(z,\theta,\theta')dz, \text{ where } \\
& \vphi(z,\ga)= \log |z+\ga|^2- \log |z-\ga|^2, \text{ and } H(z,\theta,\theta') \text{ as in \eqref{princ2}.} 
\end{split}
\eeq
 In view of \eqref{def-Z},  $\vphi(z,\ga)$ and $H(z,\theta,\theta')$  are $C^\infty$ functions and we can use stationary phase methods to analyze the asymptotic expansion of the integral as $h\downarrow 0$.  But
\[
\begin{split}
\frac14 \nabla_z \vphi(z,\ga)= \frac{z+\ga}{2|z+\ga|^2}- \frac{z-\ga}{2|z-\ga|^2}= \frac{-2\lan z,\ga\ran z+ (|z|^2+|\ga|^2)\ga}{|z-\ga|^2\, |z+\ga|^2}.
\end{split}
\]
We find that
\beq\label{bd-der}
\frac{1}{16}\left|\nabla_z \vphi(z,\ga)\right|^2=  \frac{|\ga|^2}{|z-\ga|^2\, |z+\ga|^2}. 
\eeq
We define the vector field
\[
\begin{split}
& \mcz(z,\ga)= \frac{1}{|\nabla_z \vphi|^2}\lan \nabla_z \vphi, \p_z\ran = \frac{1}{|\nabla_z \vphi|^2}\sum_{j=1}^{n+1} \p_{z_j} \vphi \p_{z_j}= |\ga|^{-1} \mcr(z,\ga), \\
& \text{ where } \mcr(z,\ga) = \left( -\lan z,\frac{\ga}{|\ga|}\ran z+ (|z|^2+|\ga|^2) \frac{\ga}{|\ga|}\right)\cdot \nabla_z \\
\end{split}
\]
and of course $\mcr$  is a  vector field with  $C^\infty$ coefficients in $z$. Then
\[
\mcz e^{\ioh \vphi(z,\ga)} =\ioh e^{\ioh \vphi(z,\ga)} \text{ and } \;\  \mcr e^{\ioh \vphi(z,\ga)} =i\frac{|\ga|}{h} e^{\ioh \vphi(z,\ga)} \text{ and } 
\]
and so for $N\in \mn$,
\[
\begin{split}
& (\ioh |\ga|)^{N} \Upsilon_V(h,\theta,\theta')= \int_{\mr^{n+1}} \left( \mcr(z,\ga)^N e^{\ioh \vphi(z,\ga)}\right) H(z,\theta,\theta') dz= \\
&  \int_{\mr^{n+1}}   e^{\ioh \vphi(z,\ga)} (\mcr(z,\ga)^T)^NH(z,\theta,\theta') dz,
\end{split}
\]
where $\mcr^T$ is the transpose of $\mcr$.  Therefore, if $|\ga|>r>0$,
\[
\left| (\ioh |\ga|)^{N} \Upsilon_V(h,\theta,\theta')\right| \leq C_N
\]
and we conclude that
\[
\text{ if } |\ga|> \del>0 , \text{ then }  \left| \Upsilon_V(h,\theta,\theta')\right| \leq C(r,N) h^{N} \text{ for all } N\in \mn.
\]
This means that $\Upsilon_V(h,\theta,\theta')$ decays rapidly in $h$ is $|\theta-\theta'|>\del>0.$

From now on we assume that $|\ga|<r$, with $r$ small enough.  This implies that $\ha \left| \theta+\theta'\right|>1-\frac{r}{2}$ and  since $V(w)=0$ for $|w|>1-\eps$, we deduce that,  on the support of $V$,
\[
|z|=\left|w-\ha(\theta+\theta')\right|> \ha\left|\theta+\theta'\right|-|w|>  \eps-\frac{r}{2}>\frac{r}{4}>0, \text{ for } r \text{ small enough.}
\]
Therefore, we write
\[
\begin{split}
& \vphi(z,\ga)= \log|z+\ga|^2-\log|z-\ga|^2 = 4\lan \frac{z}{|z|^2}, \ga\ran  - \lan \ga, \mathfrak{G}(z,\ga) \ga\ran, \text{ where } \\
&  (\mathfrak{G}(z,\ga) \ga)_j=\sum_{k}\mathfrak{G}_{jk}(z,\ga) \ga_k, \;\  \mathfrak{G}_{jk}(z,\ga) \in C^\infty \text{ and  for any } \alpha\in \mn^{n+1}, \\
&  |D_z^\alpha \mathfrak{G}(z,\ga)| \leq C_\alpha \text{ for } z\in \Omega  \text{ compact and } |z|>\frac{r}{4}
\end{split}
\]
So we write
\[
\begin{split}
& e^{\ioh \vphi(z,\ga)} - e^{-\frac{4i}{h} \lan \frac{z}{|z|^2}, \ga\ran} = e^{-\frac{4i}{h}\lan \frac{z}{|z|^2}, \ga\ran}\left(e^{\ioh \lan \ga, \mathfrak{G}(z,\ga) \ga\ran}-1\right)= \\
 &  e^{-\frac{4i}{h}\lan \frac{z}{|z|^2}, \ga\ran}  \sum_{m=1}^\infty\frac{1}{m!} \left(\ioh \lan \ga, \mathfrak{G}(z,\ga)  \ga \ran\right)^m
\end{split}
\]
So we conclude that
\[
\begin{split}
& \Upsilon_V(h,\theta,\theta')= \int_{\mr^{n+1}} e^{-\frac{4i}{h}\lan \frac{z}{|z|^2}, \ga\ran}H(z,\theta,\theta') dz + \\
& \sum_{m=1}^\infty\frac{1}{m!}  \int_{\mr^{n+1} }  e^{-\frac{4i}{h}\lan \frac{z}{|z|^2}, \ga\ran}  \left(\ioh \lan \ga, \mathfrak{G}(z,\ga) \ga \ran\right)^m  H(z,\theta,\theta') \ dz.
 \end{split}
 \]
 Now recall that on the support of the integrand $R> |z|>\frac{r}{4}$ and so we make a change of variables
 \[
 y=\frac{z}{|z|^2} \text{ and so } z=\frac{y}{|y|^2}, \text{ and } \frac{4}{r} > |y|> \frac{1}{R}>0,
 \]
 and the integral becomes
 \[
\begin{split}
& \Upsilon_V(h,\theta,\theta')= \int_{\mr^{n+1}} e^{-\frac{4i}{h}\lan y, \ga\ran}H(\frac{y}{|y|^2},\theta,\theta') |y|^{-2(n+1)} dy+ \\
&  \sum_{m=1}^\infty \frac{1}{m!} \int_{\mr^{n+1}} e^{-\frac{4i}{h}\lan y, \ga\ran}  \left(\ioh \lan \ga, \mathfrak{G}(\frac{y}{|y|^2},\ga)  \ga \ran\right)^m  H(\frac{y}{|y|^2},\theta,\theta')  \ |y|^{-2(n+1)} dy.
 \end{split}
 \]
 But notice that
 \[
 \ioh \mathfrak{G}_{jk}(\frac{y}{|y|^2},\ga) \ga_j\ga_k  e^{-\frac{4i}{h} \lan y,\ga\ran} = -\frac{ih}{4^2} \mathfrak{G}_{jk}(\frac{y}{|y|^2},\ga) \p_{y_j}\p_{y_k}e^{-\frac{4i}{h} \lan y,\ga\ran}.
 \]
 and therefore
 \[
 \begin{split}
& \left( \ioh \mathfrak{G}_{jk}(\frac{y}{|y|^2},\ga) \ga_j\ga_k \right)^m e^{-\frac{4i}{h} \lan y,\ga\ran}= \left(-i\frac{h}{4^2}\right)^m  P_m(y,\ga,D) e^{-\frac{4i}{h}\lan y,\ga\ran}, \text{ where }  \\
 & P_m(y,\ga,D) =   \sum \mathfrak{G}_{j_1k_1} \mathfrak{G}_{j_2k_2}\ldots \mathfrak{G}_{j_mk_m} \p_{y_{j_1}}\p_{y_{k_1}} \p_{y_{j_2}}\p_{y_{k_2}}
 \ldots \p_{y_{j_m}}\p_{y_{k_m}}.
 \end{split}
 \]
  and  then integration by parts gives
 \[
 \begin{split}
 & \int_{\mr^{n+1}}   e^{-\frac{4i}{h}\lan y, \ga\ran} \left(\ioh  \lan \ga, \mathfrak{G}(\frac{y}{|y|^2},\ga)  \ga \ran\right)^m  H(\frac{y}{|y|^2},\theta,\theta') \ |y|^{-2(n+1)} dy = \\
&  \left(-\frac{ih}{4^2}\right)^m  \int_{\mr^{n+1}}  H(\frac{y}{|y|^2},\ga) P_m(y,\ga,D) e^{-\frac{4i}{h}\lan y, \ga\ran}  \ |y|^{-2(n+1)} dy = \\
& \left(-\frac{ih}{4^2}\right)^m  \int_{\mr^{n+1}}   e^{-\frac{4i}{h}\lan y, \ga\ran} P_m(y,\ga,D)^T\left(H(\frac{y}{|y|^2},\ga) 
\ |y|^{-2(n+1)}\right) dy,
\end{split}
\]
where $P_m(y,\ga,D)^T$ denotes the transpose of $P_m(y,\ga,D)$.
On the other hand,  it is shown in \cite{ErdTri} that
\[
c(\novt+\frac{i}{h}) c(\novt-\frac{i}{h}) \sim \frac{1}{4\pi^n}h^{2-n}+ \sum_{j=1}^\infty d_j h^{2-n+j}
\]
and so we  have proved the Lemma.
\end{proof}

Next we further simplify the formula \eqref{princ2} and establish the connection with the geodesic $X$-ray transform of $V$.
\begin{prop}  If $\del>0$, the Schwartz kernel of  $\mcu_V$ satisfies
\beq\label{formula1}
\begin{split}
& \frac{1}{h} \mcu_V(\novt+\frac{i}{h},\theta,\theta')= O(h^\infty), \text{ provided }  |\theta-\theta'|>\del,  \text{ and otherwise }  \\
& \frac{1}{h} \mcu_V(\novt+\frac{i}{h},\theta,\theta')=   \frac{1}{(2\pi h)^n}\int_{\Pi(\theta)} e^{-\frac{i}{h}\lan \xi,\theta-\theta'\ran} \biggl(\mathfrak{H}(\theta,\xi) + h \mathfrak{R}(h,\theta,\xi) \biggl) \ d\mu_{\theta}(\xi) + O(h^\infty), \\
& \text{ where }  \Pi(\theta)=\{ \xi\in \mr^{n+1}: \lan \xi,\theta \ran=0\}  \text{ and } \mu_{\theta}  \text{ is the Lebesgue measure on the plane }  \Pi(\theta),  \\
& \mathfrak{H}(\theta,\xi)= 2^{n}\int_{\mr}  H( \frac{t\theta+\frac{\xi}{2}}{t^2+\frac{|\xi|^2}{4}},\theta,\theta) (t^2+\frac{|\xi|^2}{4})^{-n-1}\ dt= -2^{n-1} X(V)(\frac{\xi}{2},\theta) \text{ as  defined in \eqref{RTV} } \\
& \text{ and }  \mathfrak{R}(h,\theta,\xi) \sim  \sum_{j=0}^\infty h^j  \mathfrak{R}_j(\theta,\xi), \;\ \mathfrak{R}_j(\theta,\xi) \in C_0^\infty( \msn \times \mr^{n+1}), \text{ with support independent of } j\\
& \text{ and } O(h^\infty)= \int_{\Pi(\rho)} e^{-\frac{i}{h}\lan \xi,\theta-\theta'\ran} \mathfrak{a}(h, \theta,\xi) \ d\mu_{\theta'}(\xi), \text{ for some } \mathfrak{a}(h, \theta,\xi) \\
& \text{ which is compactly supported and }  |(\p_\theta,\p_\xi)^\alpha \mathfrak{a}(h, \theta,\xi)  | \leq C_{N,\alpha} h^N \text{ for all } N, \alpha.
\end{split}
\eeq
\end{prop}
\begin{proof}
As observed in Lemma \ref{princ1}, if $\del>0$ and $\chi(s)\in C_0^\infty(\mr)$, $\chi(s)=1$ if $|s|<\del$ and $\chi(s)=0$ if $|s|>2\del$, then
\[
\begin{split}
& \frac{1}{h}  \mcu_V(\novt+\frac{i}{h},\theta,\theta')= \\
  \frac{2^n}{ (\pi h)^n}\int_{\mr^{n+1}} e^{-\frac{2i}{h}\lan y, \theta-\theta'\ran}\chi(|\theta-\theta'|^2) & \ \left(H(\frac{y}{|y|^2},\theta,\theta') + hR(h,y,\theta,\theta')\right) |y|^{-2(n+1)} \ dy+ O(h^\infty),
\end{split}
\]
with $R(h,y,\theta,\theta')$ satisfying an expansion as in \eqref{princ2}. We proceed as in \cite{BirYaf1,BulPus1} and define
\[
\mathfrak{J}(\theta,\theta')= \frac{\theta+\theta'}{|\theta+\theta'|} ,
\]
 and we observe that since $\theta,\theta'\in \msn$, $\lan \mathfrak{J}(\theta,\theta'), \theta-\theta'\ran=0$ and so we can write
\[
y=  t \mathfrak{J}(\theta,\theta')+ \frac{\xi}{2}, \text{ with } \xi \text{ such that }  \lan \xi,\mathfrak{J}(\theta,\theta')\ran=0,
\]
and since $|\mathfrak{J}(\theta,\theta')|=1$, we have
\[
\begin{split}
&  \frac{1}{h} \mcu_V(\novt+\frac{i}{h},\theta,\theta')  \ =  \frac{1}{(\pi h)^n}\int_{\Pi(\mathfrak{J})} e^{-\frac{i}{ h}\lan \xi,\theta-\theta' \ran}\biggl( \mch(\theta,\theta',\xi) +h \mcr(h,\theta,\theta',\xi)\biggr) \ d\mu_{\mathfrak{J}}(\xi) + O(h^\infty), \\
& \text{ with }  \Pi(\mathfrak{J})  \text{ and } \mu_{\mathfrak{J}} \text{ as defined in \eqref{formula1}  with } \mathfrak{J} \text{ in place of } \theta, \\ 
& \mch(\theta,\theta',\xi)= \int_{\mr} \chi(|\theta-\theta'|^2) H( \frac{t\mathfrak{J}(\theta,\theta')+\frac{\xi}{2}}{t^2+\frac{|\xi|^2}{4}},\theta,\theta') (t^2+\frac{|\xi|^2}{4})^{-n-1}\ dt, 
\text{ and } \\
& \mcr(\theta,\theta',\xi)= \int_{\mr} \chi(|\theta-\theta'|^2) R( h, t\mathfrak{J}(\theta,\theta') + \frac{\xi}{2},\theta,\theta')\ dt, 
\end{split}
\]


As in the proof of Lemma \ref{princ1}, we write the Taylor series expansion  of $\mch(\theta,\theta',\xi)$  and $\mcr_j(\theta,\theta',\xi)$  along the diagonal $\{\theta=\theta'\}$ 
\[
\begin{split}
& \mch(\theta,\theta',\xi)\sim \mch(\theta,\theta,\xi)+ \sum_{|\alpha|=1} (\theta-\theta')^\alpha T_\alpha(h,\theta,\xi), \text{ with } T_\alpha \in C^\infty, \\
&\mcr_j(\theta,\theta',\xi) \sim \mcr_j(\theta,\theta,\xi)+ \sum_{|\alpha|=1}^\infty (\theta-\theta)^\alpha \mct_{\alpha,j}(h,\theta,\xi), \text{ with } \mct_{\alpha,j} \in C^\infty, \\
\end{split}
\]
 use that
\[
(\theta-\theta')e^{-\frac{i}{h} \lan \theta-\theta',\xi\ran}= -\frac{h}{i} \nabla_\xi e^{-\frac{i}{h} \lan \theta-\theta', \xi\ran},
\]
and that $\chi(|\theta-\theta'|^2)=1$ near $\{\theta=\theta'\}$ and integrate by parts in $\xi$ and use that $\mathfrak{J}(\theta,\theta)=\theta$  to conclude that
\[
\begin{split}
&  \frac{1}{h} \mcu_V(\novt+\frac{i}{h},\theta,\theta')=  \frac{1}{(\pi h)^n}\int_{\Pi(\mathfrak{J})} e^{-\frac{i}{ h}\lan \xi,\theta-\theta' \ran}\biggl( \mathfrak{H}(\theta,\xi) +h \mathfrak{R}(h,\theta,\xi)\biggr) \ d\mu_{\mathfrak{J}}(\xi) + O(h^\infty), \\
&  \mathfrak{H}(\theta,\xi)= \mch(\theta,\theta,\xi), \;\  \mathfrak{R}(h,\theta) \sim \sum_{j=0}^\infty h^j \mathfrak{R}_j(\theta,\xi),\;\  \mathfrak{R}_j \in C_0^\infty, \text{ with support independent of } j
\end{split}
\]

This is not quite the formula we want because the integral is over the hyperplane $\Pi(\mathfrak{J})$ and we want it over $\Pi(\theta)$ and this is clarified by the following 
\begin{lemma}  Let $\theta,\ga\in \msn$ be such that $|\theta-\ga| <\eps$ is small enough.  Let $\sigma=\frac{\theta+\ga}{|\theta+\ga|}$, $\Pi(\bullet)$ and $\mu_\bullet$ be defined as in \ref{formula1}, $\bullet=\theta, \sigma$. If $F(z)\in C_0^\infty(\mr^{n+1})$, then
 \[
 \begin{split}
&  \int_{\Pi(\sigma)} e^{\ioh \lan z,  \ga-\theta\ran} F(z) \ d\mu_{\sigma}(z)-  \int_{\Pi(\theta)} e^{\ioh \lan z, \ga-\theta\ran} F(z) \ d\mu_{\theta}(z) = \\
& h   \int_{\Pi(\theta)} e^{\ioh \lan z,  \ga- \theta\ran} \mathfrak{F} (h,z,\theta) \ d\mu_{\theta}(z) + O(h^\infty),  \text{ where } \\
&  \mathfrak{F}(h,z, \theta)\sim \sum_{j=0}^\infty h^j \mathfrak{F}_j(z,\theta), \;\ \mathfrak{F}_j \in C_0^\infty, \text{ with support independent of } j
 \end{split}
 \]
\end{lemma}
\begin{proof}   If  $A$ is a linear orthogonal transformation and $\ga=A\ga^*$ and $\theta= A \theta^*$, we have
\[
\begin{split}
& \int_{\Pi(\sigma)} e^{\ioh \lan z,  \ga-\theta\ran} F(z) \ d\mu_{\sigma}(z)-  \int_{\Pi(\theta)} e^{\ioh \lan z, \ga-\theta\ran} F(z) \ d\mu_{\theta}(z) = \\
 & \int_{\Pi(\sigma^*)} e^{\ioh \lan z,  \ga^*-\theta^*\ran} F(Az) \ d\mu_{\sigma^*}(z)-  \int_{\Pi(\theta^*)} e^{\ioh \lan z, \ga^*-\theta^*\ran} F(Az) \ d\mu_{\theta^*}(z),
\end{split}
\]
 and so we may assume that $\theta=(0,1)$ and $\ga= (\ga', \ga_{n+1})$, with $|\ga'|<\eps$ and $|1-\ga_{n+1}|<\eps$, for $\eps$ small enough. If we denote $z=(z',z_{n+1})$, $z'=(z_1, \ldots z_n)$, and since in this case $\sigma=\frac{1}{|(\ga',1+\ga_{n+1})|} (\ga',1+\ga_{n+1}),$ we have
\[
\Pi(\theta)= \{z \in \mr^{n+1}: z_{n+1}=0\} \text{ and } \Pi(\sigma)=\{ z \in \mr^{n+1}: z_{n+1}= - \frac{1}{1+\ga_{n+1}} \lan z',\ga'\ran,
\]
and so
\[
\begin{split}
&  \int_{\Pi(\theta)} e^{\ioh \lan z, \ga-\theta\ran} F(z) \ d\mu_{\theta}=  \int_{\mrn} e^{\ioh \lan z',  \ga' \ran} F(z',0) \ d z' \text{ and } \\
&  \int_{\Pi(\sigma)} e^{\ioh \lan z, \ga-\theta\ran} F(z) \ d\mu_{\sigma}=   \int_{\mrn} e^{\ioh \lan z',  \ga' \ran}  e^{-\ioh\bigl( \frac{\ga_{{}_{n+1}-1}}{\ga_{{}_{n+1}+1}} \lan z',\ga'\ran\bigr)} G(z',\ga) \ dz'= \\
&  \int_{\mrn}   e^{-\ioh \bigl(\frac{2}{\ga_{{}_{n+1}+1}}\bigr) \lan z',\ga'\ran} G(z', \ga) \ dz' \\
& \text{ where } G(z',\ga)=  F\bigl(z', - \frac{1}{1+\ga_{n+1}} \lan z',\ga'\ran\bigr)\biggl(1+ \bigl(\frac{\lan z',\ga'\ran}{1+\ga_{n+1}}\bigr)^2\biggr)^\ha, \\
 \end{split}
 \]
If we set $z'= \ha(\ga_{n+1}+1)w'$, then
\[
\begin{split}
&  \int_{\Pi(\sigma)} e^{\ioh \lan z, \ga-\theta\ran} F(z) \ d\mu_{\sigma}- \int_{\Pi(\theta)} e^{\ioh \lan z, \ga-\theta\ran} F(z) \ d\mu_{\theta}=  \\
&  \int_{\mrn} e^{\ioh \lan w',\ga'\ran }\biggl(  \bigl(\ha(\ga_{n+1}+1)\bigr)^n G\bigl( \ha(\ga_{n+1}+1)w',\ga \bigr) - F(w',0) \biggr) \ dw'= \\
& \int_{\mrn} e^{\ioh \lan w',\ga'\ran }\biggl(  \bigl(\ha(\ga_{n+1}+1)\bigr)^n  \bigl(1+ \bigl(\ha \lan \omega',\ga'\ran \bigr)^2\bigr)^\ha F\bigl(\ha(\ga_{n+1}+1)w',-\ha \lan w',\ga'\ran\bigr) - F(w',0) \biggr) \ dw'.
\end{split}
\]

Now we observe that
\[
\begin{split}
&\bigl(1+ \bigl(\ha \lan \omega',\ga'\ran \bigr)^2\bigr)^\ha= 1+\sum_{j=1}^\infty C_j \bigl(\lan \omega',\ga'\ran \bigr)^{2j}, \text{ and } \\
&  F\bigl(\ha (\ga_{n+1}+1) w',-\ha \lan w',\ga'\ran\bigr) \sim   F\bigl(\ha(\ga_{n+1}+1) w',0 \bigr)  +  \\
&  \sum_{j=1}^\infty \bigl(\p_{w_{n+1}}^j  F\bigr) \bigl(\ha(\ga_{n+1}+1) w',0\bigr) \bigl(-\ha \lan w',\ga'\ran\bigr)^j,
\end{split}
\]
and so
\[
\begin{split}
 \bigl(1+ \bigl(\ha(\lan \omega',\ga'\ran)\bigr)^2\bigr)^\ha & \, F\bigl(\ha(\ga_{n+1}+1) w',-\ha \lan w',\ga'\ran\bigr) \sim  F\bigl(\ha(\ga_{n+1}+1) w',0 \bigr) + \\
 & \sum_{j=1}^\infty F_j\bigl(\ha(\ga_{n+1}+1)w') \bigl(\lan w',\ga'\ran\bigr)^j.
\end{split}
\]

We then recognize that
\[
\bigl( \lan w',\ga'\ran\bigr)^j e^{\ioh \lan w',\ga'\ran}= \bigl(h \lan w',D_{w'} \ran \bigr)^j e^{\ioh\lan w',\ga'\ran},
\]
and therefore, after integrating by parts we obtain
\[
\begin{split}
& \int_{\mrn} e^{\ioh \lan w',\ga'\ran}  \bigl(\ha(\ga_{n+1}+1)\bigr)^n   F_j\bigl(\ha(\ga_{n+1}+1) w'\bigr) \bigl( \lan w',\ga'\ran\bigr)^j \ dw'= \\
& h^j  \int_{\mrn} e^{\ioh \lan w',\ga'\ran}   \bigl(\ha(\ga_{n+1}+1)\bigr)^n  \bigl(h \lan w',D_{w'} \ran^T \bigr)^j \bigl(  F_j \bigl(\frac{\ga_{n+1}+1}{2} w')\bigr) \ dw',
\end{split}
\]
and here $T$ denotes the transpose of the differential operator.    Since  $ \ha(\ga_{n+1}+1) = 1+ \ha( \ga_{n+1}-1)$, we may write
\[
\begin{split}
& \bigl( \lan w', D_{w'}\ran^T\bigr)^j  F_j \bigl(\ha(\ga_{n+1}+1) w')\bigr) \defn \mch_j\bigl( \ha( \ga_{n+1}+1) w'\bigr)= \mch_j\bigl( w'+ \ha( \ga_{n+1}-1) w'\bigr) \sim  \\ & \mch_j(w') + \sum_{|\alpha|=1}^\infty \mch_{j,\alpha}(w')  (\ga_{n+1}-1)^{|\alpha|} {w'}^\alpha.
\end{split}
\]
But

\[
\biggl(\frac{\ga_{n+1}+1}{2}\biggr)^n = 1+ \sum_{k=1}^n \binom{n}{k}  \biggl(\ha( \ga_{n+1}-1)\biggr)^k,
\]
and so
\[
 \bigl(\ha(\ga_{n+1}+1)\bigr)^n  \bigl(h \lan w',D_{w'} \ran^T \bigr)^j \bigl(  F_j \bigl(\ha(\ga_{n+1}+1) w')\bigr) \sim \mch_j(w')+ \sum_{j+|\alpha|=1}^\infty \wt  \mch_{j,\alpha}(w') (\ga_{n+1}-1)^{|\alpha|+j} {w'}^\alpha.
 \]

But since $\ga\in \msn$, we have
\[
\begin{split}
& \ga_{n+1}-1= \sqrt{1-|\ga'|^2}-1= \sum_{m=1}^\infty C_m |\ga'|^{2m}, \text{ and so } \\
& \bigl(\ga_{n+1}-1\bigr)^j = \sum_{m=j}^\infty C_{j,m} |\ga'|^{2m},
\end{split}
\] 
and hence
\[
 \bigl(\ha(\ga_{n+1}+1)\bigr)^n  \bigl(h \lan w',D_{w'} \ran^T \bigr)^j \bigl(  F_j \bigl(\ha(\ga_{n+1}+1) w')\bigr) \sim \mch_j(w')+ \sum_{j+|\alpha|=1}^\infty \wt  \mca_{j,\alpha}(w') |\ga'|^{2(|\alpha|+j)} {w'}^\alpha.
 \]

Finally notice that, since 
\[
|\ga'|^{2m} e^{\ioh \lan w',\ga'\ran} = (h^2 \Delta_{w'})^m  e^{\ioh \lan w',\ga'\ran} ,
\]

after integrating by parts, we conclude that
\[
\begin{split}
& \int_{\mrn} e^{\ioh \lan w',\ga'\ran}   \bigl(\ha(\ga_{n+1}+1)\bigr)^n  \bigl(h \lan w',D_{w'} \ran^T \bigr)^j \bigl(  F_j \bigl(\ha(\ga_{n+1}+1) w')\bigr) \ dw' =  \\
& \int_{\mrn}  e^{\lan w', \ga'\ran } \mathfrak{A}_j(h, w') \ dw' +O(h^\infty),  \text{ where }  \mathfrak{A}_j(h, w') \sim \mch_j(w') + \sum_{k=1}^\infty \mathfrak{A}_{jk}(w') h^k.
\end{split}
\]

 We are left to estimate the difference
\[
\int_{\mrn} e^{\ioh \lan w',\ga'\ran} \biggl( \bigl(\frac{\ga_{n+1}+1}{2}\bigr)^n  F\bigl(\frac{\ga_{n+1}+1}{2} w',0 \bigr) - F(w',0)\biggr) \ dw',
\]
which of course can be handled in the same way. 

This ends the proof of the Lemma and it also ends the proof of the Proposition.
\end{proof}
\end{proof}  

We define a very special class of semiclassical pseudodifferential operators and refer the reader to \cite{Mar,Zwo}  for  much deeper discussions about such operators.  We say that
\[
A: C_0^\infty(\mrn) \longrightarrow C_0^\infty(\mrn)
\]
is  a semiclassical pseudodifferential operator if
\beq
\label{def-sc}  A(z,hD) u(z)=  \frac{1}{(2\pi h)^n} \int_{\mrn} \int_{\mrn} e^{\ioh \lan z-y,\xi\ran} a(h,z,y,\xi)  u(y) \ dy d\xi, 
\eeq
such that
\[
\begin{split}
& a(h,z,y,\xi) \sim \sum_{j=0}^\infty h^j a_j(z,y,\xi) \text{ where } a_j(z,y,\xi)  \\ 
& \text{ is supported on a compact set which is independent of } j. 
\end{split}
\]
 In what follows we will say that
\[
 A(z,hD) =  A_0(z,hD)+ O(h^k),
 \]
 if  there exists a semiclassical pseudodifferential operator  $ B(x,hD)$  such that 
 \[
    A(z,hD)= A_0(z,hD)+ h^k B(z,hD), 
\]

We use $O(h^\infty)$ when this is true for every $k$.
Using that
\[
\begin{split}
& a_j(z,y,\xi)\sim \sum_{\alpha\in \mn} \frac{1}{\alpha!} \p_x^\alpha a_j(z,z,\xi) (z-y)^\alpha \text{ and } \\
& (z-y)^\alpha e^{\ioh \lan z-y,\xi\ran}= (h D_\xi)^\alpha e^{\ioh \lan z-y,\xi\ran},
\end{split}
\]
one can use and Borel summation formula to show that if $A(z,hD)$ is defined as in \eqref{def-sc}, then 
\[
\begin{split}
   & A(z,hD) u(z)= \frac{1}{(2\pi h)^n}  \ \int_{\mrn} \int_{\mrn} e^{\ioh \lan z-y,\xi\ran} \wt a (h,z,\xi)  u(y)  \ dy d\xi + O(h^\infty), \text{ where } \\
 & \wt a(h,z,\xi)\sim \sum_{j=0}^\infty h^j \wt a_j(z,\xi), \text{ for some }  \wt a_j(z,\xi) \in C_0^\infty,  \text{ with support independent of } j \\
 & \text{ and } \wt a_0(z,\xi)= a_0(z,z,\xi)  \text{ is defined to be the (left) principal symbol of } A(x,hD).
\end{split}
\]
  
  We can define semiclassical pseudodifferential operators on a manifold by using local coordinates, see for example \cite{Zwo}.  If $M$ is a $C^\infty$ manifold, $\{U_\mu\}$, $U_\mu\subset M$ form an open cover of $M$ and   $\Psi_\mu: U_\mu \longmapsto \mrn$,  are local charts and $\{\chi_\mu\}$ is a partition of unity subordinated to $\{U_\mu\}$, then $P$ is a semiclassical pseudodifferential operator on $M$ if for $f\in C^\infty(M)$,
  \beq\label{pdo-sym}
  \begin{split}
&  \chi_\mu P \chi_\nu f= O(h^\infty), \;\ \text{ provided } \supp \chi_\mu \cap \supp \chi_\nu=\emptyset, \text{ otherwise } \\
 & (\chi_\mu P \chi_\nu f)(\Psi_\mu^{-1}(z))=  \frac{1}{(2\pi h)^n} \int_{\mrn} e^{-\ioh \lan z-y,\xi\ran}  \chi_\mu(\Psi_\mu^{-1}(z)) a_{\mu,\nu}(h,z,\xi) (\chi_\nu f)(\Psi_\mu^{-1}(y)) \ dy d\xi,\\
& \text{ with }  a_{\mu,\nu}(h,z,\xi) \sim \sum_{j=0}^\infty h^j a_{j, \mu,\nu}(z,\xi), \text{ with } a_{j,\mu,\nu} \text{ with compact supported  independent of } j  \\
 \end{split}
  \eeq
and this characterization is independent of the atlas  $\{(U_\mu, \Psi_\mu)\}$, see for example Section 14.2 of \cite{Zwo}.

We want to show that $\ooh\mcu_V(\novt+\ioh)$ fits this definition.  A similar discussion for regular pseudodifferential operators is carried out in Section 12 of Chapter 0 of \cite{Yaf}.  We already know it satisfies the first part of \eqref{pdo-sym}, and we need to check it satisfies the second part.  For fixed $\omega\in \msn$,  let  $\Pi(\omega)$ be the plane orthogonal to $\omega$ passing through the origin:
\[
\Pi(\omega)= \{ \xi\in \mr^{n+1}: \lan \xi,\omega\ran=0\},
\]
 and let
\[
 \Omega_{\omega}= \{\theta \in \msn: |\theta-\omega|<\ha\}
 \]
 and let $\Psi_\omega:  \Omega_{\omega}\longmapsto U_\omega $ denote the orthogonal projection of $\msn$ onto $\Pi(\omega)$ restricted to $\Omega_{\omega}$. For each fixed $\omega,$ its inverse   is given by
   \[
   \begin{split}
   &  \Psi_\omega^{-1}: U_{\omega}  \longmapsto \Omega_{\omega} \\
   & \Psi_\omega^{-1}(y)= y+(1-|y|^2)^\ha \omega.
   \end{split}
   \]
   This collection $\bigl\{\bigl(\Omega_{\omega}, \Psi_\omega\bigr), \omega\in \msn \bigr\}$ gives an atlas of $\msn$.    Of course one just needs finitely many $\omega\in \msn$ to produce form an atlas.

    In our case we have an operator of the form 
  \[
  \begin{split}
  & \bigl(A(\theta, h D) f\bigr)(\theta)= \frac{1}{(2\pi h)^n} \int_{\msn} \int_{\Pi(\theta)} e^{\ioh \lan \xi, \theta-\theta'\ran } a(h, \theta,\xi) f(\theta') \ d\mu_\theta(\xi) d\theta', \\
  & \text{ with } a(h, \theta,\xi)\sim \sum_{j=0}^\infty h^j a_j (\theta,\xi),  \;\ a_j \in C_0^\infty,
  \end{split}
  \]
   and we want to show that for each fixed $\omega\in \msn$,  $(\chi A(\theta, hD) \chi f)(\Psi_\omega^{-1})$ is given by a formula as in \eqref{pdo-sym}. In fact we just need to verify that this is the case if $\theta=\omega$.  Indeed,  if  $T=T_{\omega,\theta}$ is the orthogonal transformation such that $\theta=T\omega$, and if we define  $\theta'=T\omega'$, one can check that 
 \[
\bigl(A(T\omega, hD) f )(T\omega)= \frac{1}{(2\pi h)^n} \int_{\msn} \int_{\Pi(\omega)} e^{\ioh \lan \xi, w-w'\ran} a(h, T\omega, T\xi) f(T w') \ d\mu_\omega(\xi) d\omega'.
\]
So we just need to consider the case $\omega=\theta$ and verify that $(\chi A(\theta, hD) \chi f)(\Psi_\theta^{-1})$  is given by a formula as in \eqref{pdo-sym}.  We set $z= \Psi_\theta(\theta)$  (we are just shifting the origin in $\Pi(\theta)$ to $z$ instead of $0$)  and   
  $y= \Psi_\theta(\theta'),$ and $\chi\in C^\infty(\msn)$ supported near $\theta,$ then $\theta'=\Psi_\theta^{-1}(y)=y+(1-|z-y|^2)^\ha \theta$, and since $\lan \xi,\theta\ran=0$ for $\xi\in \Pi(\theta),$  
  \[
  \begin{split}
  & (\chi A(\theta, h D) \chi f)(\Psi_\theta^{-1}(z))= \\
  \frac{1}{(2\pi h)^n}\int_{\Pi_\theta}  \int_{\Pi_\theta} & \ e^{\ioh \lan z-y, \xi\ran } \chi(\Psi_\theta^{-1}(z)) a(h, \Psi_\theta^{-1}(z), \xi) (\chi f)(\Psi_{\theta}^{-1}(y)) \frac{1}{\sqrt{1-|y-z|^2}} \ d\mu_\theta(\xi) d\mu_\theta(y).
  \end{split}
  \]
  Here we used that $\frac{1}{\sqrt{1-|y-z|^2}} dy = d\theta'.$   If we proceed as above we can see that the integral decays is of order $O(h^\infty)$ if $|z-y|>\del>0,$ and so we may  expand 
  \[
  \frac{1}{\sqrt{1-|y-z|^2}}=1+ \sum_{j=0}^\infty C_j |z-y|^{2j}, \text{ for } |z-y|<\eps,
  \]
  and  use that $|z-y|^2 e^{\ioh \lan z-y,\xi\ran}= h^2 \Delta_\xi e^{\ioh \lan z-y,\xi\ran}$ and integrate by parts in $\xi$, we find that
  \[
\begin{split}  
& (\chi A(\theta, h D) \chi f)(\Psi_\theta^{-1}(z))= \frac{1}{(2\pi h)^n}\int_{\mrn}  \int_{\mrn} e^{\ioh \lan z-y, \xi\ran } b(h, z, \xi) (\chi f)(\Psi_{\theta}^{-1}(y)) \ d\xi dy + O(h^\infty),  \\
& \text{ where } b(h,z,\xi)\sim  a_0(\Psi_\theta^{-1}(z), \xi) + \sum_{j=1}^\infty h^j b_j(z,\xi), \;\ b_j\in C_0^\infty,
\end{split}
  \] 
and we are just using $d\xi=d\mu_\theta(\xi)$ and $d y=d\mu_\theta(y)$.  In our case $a_0(\Psi_{\theta}^{-1}(z),\xi)=\mathfrak{H}(\theta,\xi)$, defined in \eqref{formula1} and \eqref{RTV}.  Notice that any other operator $\wt A$ of the same type such that  
 \[
  \wt b(h,\Psi_\theta^{-1}(z),\xi)\sim  \mathfrak{H}(\theta,\xi) + \sum_{j=1}^\infty h^j \wt b_j(\theta,\xi), \;\ \wt b_j\in C_0^\infty,
  \]
  will satisfy $A-\wt A= O(h),$ and so we say that $ \mathfrak{H}(\theta,\xi)$ is the principal symbol of $A(\theta,hD)$, and  this proves the  following 
 
\begin{theorem}\label{princ4}  If $V\in C_0^\infty(\mb^{n+1})$,  the operator $\frac{1}{h} \mcu_V(\novt+\frac{i}{h})$ defined in \eqref{mce-st} is a semiclassical pseudodifferential operator on $\ms^n$ and its principal symbol is equal to $\mathfrak{H}(\theta,\xi)=-2^{n-1} X(V)(\frac{\xi}{2},\theta)$. 
\end{theorem}

\section{The Proof of Theorem \ref{main} for $V\in C_0^\infty(\mhn)$} \label{TRP}

We first prove Theorem \ref{main} for $V\in C_0^\infty(\mhn)$ and for polynomials $f(t)$ and then use the Stone-Weiestrass theorem to extend it to compactly supported functions.  As already mentioned in Section \ref{OUT}, the idea of using the Stone-Weiestrass theorem  in this context is not new.  In the next section  we extend the result for $V\in x^m L^\infty(\mhn)$, $m>1$.

If $V\in C_0^\infty(\mhn)$, we have shown that $\frac{1}{h} \mcu_V(\novt+\frac{i}{h})$ is a semiclassical pseudodifferential operator whose principal symbol is $-2^{n-1} X(V)(\frac{\xi}{2},\theta)$.  The calculus of semiclassical pseudodifferential operator, see \cite{Mar,Zwo} and Theorem \ref{princ4}  give that for any $p\in \mn$, $(\frac{1}{h} \mcu_V(h))^p$ satisfies
\[
\begin{split}
& \left( \frac{1}{h} \mcu_V(\novt+\frac{i}{h})\right)^p(\theta,\theta')= \frac{1}{(2\pi h)^n}\int_{\Pi(\theta)} e^{-\frac{i}{h}\lan \xi,\theta-\theta' \ran} \left(-2^{n-1} X(V)(\frac{\xi}{2},\theta)\right)^p d\mu_\theta(\xi) +  \\
&  \frac{h}{(2\pi h)^n}\int_{\Pi(\theta)} e^{-\frac{i}{h}\lan \xi,\theta-\theta' \ran} b(h,\theta,\xi) \ d\xi + O(h^\infty), \;\   b(h,\theta,\xi) \sim \sum_{j=0}^\infty h^j b_j(\theta,\xi).
\end{split}
\]
and therefore its trace satisfies
\beq\label{TR0}
\begin{split}
& (2\pi h)^n Tr\left( \frac{1}{h} \mcu_V(\novt+\frac{i}{h})\right)^p= \int_{\msn}(2\pi h)^n \left( \frac{1}{h} \mcu_V(\novt+\frac{i}{h}h)\right)^p(\theta,\theta) \ d\theta+ O(h)=\\
&  \int_{\msn} \int_{\Pi(\theta)} \left(-2^{n-1} X(V)(\frac{\xi}{2},\theta)\right)^p \  d\mu_{\theta}(\xi) d\theta + O(h)= \int_{\mr} t^p \ d\nu(t) + O(h).
\end{split}
\eeq
where $\nu$ is the measure defined in \eqref{meas3}. 

Since $ V\in x^m L^\infty(\mhn)$ for all $m$, it follows from \eqref{ptr1N} and \eqref{TR0}  that for $p\geq 1$,
\beq\label{th-sc}
(2\pi h)^n \sum_{e^{i\del(s)}\in \spec \Tau_V} \bigl(\frac{1}{h} \ka_j(h)\bigr)^p - \int_{\mr} t^p d\nu(t) =O(h), \;\ \del_j(\novt+\frac{i}{h})=\ka_j(h).
\eeq

This proves Theorem \ref{main} if $V\in C_0^\infty(\mhn)$, $f(t)=t^p$, with $p\in \mn$, $p\geq 1$.   To prove the result for $f$ such that $\frac{f(t)}{t}\in C_0(\mr)$, we appeal to the Stone-Weistrass theorem.  Since $V\in C_0^\infty(\mhn)$, $X(V)(\frac{\xi}{2},\theta)=0$ for $|\xi|$ large and so the measure $\nu$ is supported on $[-R,R]$ for some $R>0$. Similarly, in view of \eqref{est-EV}, $\frac{1}{h} \ka_j(h)\in [-R,R]$, for $h$ small, for all $j$. So we may work with $f$ supported in $[-R,R]$.  Define the space
\[
\begin{split}
& \Gamma([-R,R])= \bigl\{ f \in C^0([-R,R]): \frac{f(t)}{t} \text{ is continuous } \bigr\}, \\
& \text{ equipped with the  norm } ||f||_{{}_\Gamma}= \sup_{x\in [-R,R]} \biggl| \frac{f(t)}{t} \biggr|.
\end{split}
\]
Since $V$ is compactly supported, we can take $p=1$ in \eqref{eig-ap0}, and so if $\mu_h$ is the measure defined in \eqref{count2}, then 
\[
\begin{split}
& \bigl| \lan \mu_h, f\ran\bigr|= \biggl| (2\pi h)^n \sum_{e^{i\ka_j(h)} \in \spec(\Tau_V)} f\bigl( \frac{\ka_j(h)}{h}\bigr) \bigr| \leq \\
&  (2\pi h)^n ||f||_{{}_\Gamma}\sum_{e^{i\ka_j(h)} \in \spec(\Tau_V)} | \frac{\ka_j(h)}{h} \bigr| \leq C ||f||_{{}_\Gamma} ||\frac{1}{h}\mcu(\novt+\frac{i}{h})||_{\Sigma_1}.
\end{split}
\]

Since $\frac{f(t)}{t}$ is continuous, for any $\eps>0$, there exists a polynomial $q(t)$ such that $\bigl| \frac{f(t)}{t}-q(t)\bigr|<\eps$, for all $t$, and so
$\bigl| \frac{f(t)-tq(t)}{t}\bigr|<\eps$ for all t, and we conclude that the space of polynomials $q(t)$ such that $q(0)=0$ is dense in $\Gamma([-R,R])$. But the argument used above implies that if $q$ is a such a polynomial, then 
\[
\bigl| \lan \mu_h, f-q\ran \bigr| \leq C || f-q||_{{}_\Gamma}.
\]
This implies that
\[
\biggl|\lan \mu_h, f\ran -\int_{\mr} q(t) \ d\nu(t) -O(h)\biggr|\leq C || f-q||_{{}_\Gamma}.
\]
If we first take the limit as $h\rightarrow 0$ and then take a sequence of polynomials $q_j$,  which converge to $f$ in $\Gamma([-R,R])$ and the result follows.

\section{The Proof of the General Case of Theorem \ref{main}}

 We have proved Theorem \ref{main} for $V\in C_0^\infty(\mhn)$ and we want to extend it to include the cases where  $V\in x^m L^\infty(\mhn)$, for some $m\in (1,\infty)$.  We know from \eqref{mcu} that if $V\in x^m L^\infty(\mhn)$ and $p\in \mn$ is such that   $p\geq 1$  and $p>\frac{2(n-\ha)}{m-1}$, $\bigl( \ooh\mcu_V(\novt+\frac{1}{h})\bigr)^p$ is trace class and 
\[
\biggl|Tr\biggl[ (2\pi h)^{n} \biggl( \frac{1}{h} \mcu_V(\novt+\frac{i}{h})\biggr)^p\biggr] \biggr| \leq C||x^{-m} V||_{L^\infty}.
\]
We also know from \eqref{amb1} and \eqref{amb2} that if $V_j\in x^m L^\infty(\mhn)$, and $A_j=\frac{1}{h} \mcu_{V_j}(\novt+\ioh)$, $j=1,2$, that
\[
\biggl|Tr\bigl( A_1^p-A_2^p\bigr)\biggr|\leq  \sum_{j=0}^{p-1} ||A_1||_{\Sigma_p}^j||A_1-A_2||_{\Sigma_p}|| A_2||_{\Sigma_p}^{p-j-1}
\]
and so it follows from \eqref{mcu} that
\beq\label{TR1}
h^n \biggl|Tr\bigl( A_1^p-A_2^p\bigr)\biggr|\leq  C||x^{-m} (V_1-V_2)||_{{}_{L^\infty}} \sum_{j=0}^p  ||x^{-m} V_1||_{{}_{L^\infty}}^j ||x^{-m} V_2||_{{}_{L^\infty}}^{p-j-1}.
\eeq

On the other hand,
\[
\begin{split}
& \int_{\msn} \int_{\Pi(\theta)}\biggl(X(V_1)(\frac{\xi}{2},\theta)\biggr)^p -\biggl(X(V_2)(\frac{\xi}{2},\theta)\biggr)^p \biggr)\  d_{\mu_\theta}(\xi) d\theta= \\
& \int_{\msn} \int_{\Pi(\theta)}\biggl( \biggl(X(V_1)(\frac{\xi}{2},\theta)-X(V_2)(\frac{\xi}{2},\theta)\biggr) Q(X(V_1)(\frac{\xi}{2},\theta),X(V_2)(\frac{\xi}{2},\theta))  d_{\mu_\theta}(\xi) d\theta,
\end{split}
\]
where $Q(x,y)$ is a homogeneous polynomial of degree $p-1$. Therefore, in view of \eqref{norm1} and \eqref{bdXR}, we obtain
\beq\label{TR2}
\begin{split}
& \biggl|  \int_{\msn} \int_{\Pi(\theta)}\biggl(X(V_1)(\frac{\xi}{2},\theta)\biggr)^p -\biggl(X(V_2)(\frac{\xi}{2},\theta)\biggr)^p \biggr)\   d_{\mu_\theta}(\xi) d\theta\biggr|\leq \\
& CQ( ||x^{-m}V_1||_{{}_{L^\infty}}, ||x^{-m}V_2||_{{}_{L^\infty}}) ||x^{-m}(V_1-V_2)||_{{}_{L^\infty}} \int_{\mrn} (1+|\xi|)^{-p m} \  d_{\mu_\theta}(\xi) \leq  \\
& CQ( ||x^{-m}V_1||_{{}_{L^\infty}}, ||x^{-m}V_2||_{{}_{L^\infty}}) ||x^{-m}(V_1-V_2)||_{{}_{L^\infty}}.
\end{split}
\eeq
Therefore, \eqref{TR1} and \eqref{TR2} imply that \eqref{TR0} also holds for $V\in Z_{0,m}\defn \overline{C_0^\infty(\mhn)}$, which is defined to be the closure of $C_0^\infty(\mhn)$ with the topology of $x^m L^\infty$, provided $p> \frac{2(n-\ha)}{m-1}$.

One has that
\[
x^m L^\infty(\mhn) \subset  x^{m'} L^\infty(\mhn), \text{  provided } m\geq m',
\]
but in fact one also has 
\[
x^m L^\infty(\mhn) \subset Z_{0,m'}, \text{ provided } m>m'.
\]
To see that, notice that if $V\in x^m L^\infty(\mhn)$ and $\chi_k(w)=1$ if $|w|<1-\frac{1}{k}$ and $\chi_k(w)=0$ if $|w|>1-\frac{1}{2k}$, then
\[
\begin{split}
 ||x^{m'}(V-\chi_k V)||= & ||x^{-m} V x^{m-m'}(1- \chi_k)||_{}\leq  ||x^{-m}V||_{{}_{L^\infty}} ||x^{m-m'}(1- \chi_k(w))||_{{}_{L^\infty}} \leq \\
& ||x^{-m}V||_{{}_{L^\infty}} k^{m'-m} \rightarrow 0 \text{ as } k\rightarrow \infty.
\end{split}
\]
This means that $\chi_k V\rightarrow V$ in $x^{m'}L^\infty(\mhn)$, if $m'<m$. So, if $p> \frac{2(n-\ha)}{m-1}$, pick $m'\in (1,m)$ such that 
$p> \frac{2(n-\ha)}{m'-1}$, and so the result holds for $Z_{0,m'} \supset x^m L^\infty(\mhn)$. This implies that Theorem \ref{main} holds for polynomials of degree $p> \frac{2(n-\ha)}{m-1}$.  We modify  the proof of the case $V\in C_0^\infty(\mhn)$, and we use space $\Gamma_p([-R,R])$,  where $\frac{f(t)}{t^p}$ is continuous, the argument goes through and we prove \eqref{lim-high2} for continuous functions $f$ which vanish to order $p$ at zero.

\section{The Proof of Theorem \ref{inverse}}

If $V(w)= V(\rho(w))$, where $\rho(w)= d_{g_0}(w,0)=\log\bigl(\frac{1+|w|}{1-|w|}\bigr)$,  it follows from \eqref{RTV} that
\beq\label{RTVR}
\begin{split}
 & -2^{n-1} X_V(t,\xi) = G_V(r)= \int_{-\infty}^{-\ha} V\big(T(r,t)\bigr) \frac{1}{1+2t} \ dt, \text{ where } \\
& r=|\xi|, \;\ T(r,t)=  \log\biggl(\frac{1+A(t,r)}{1-A(t,r)}\biggr) \text{ and }  A(t,r)=\biggl[ \frac{(1+t)^2+r^2}{t^2+r^2}\biggr]^\ha.
\end{split}
\eeq

Notice that $x=\frac{1-|w|}{1+|w|}=e^{-d_{g_0}(w,0)}$, defined in \eqref{g0x}, and so so if $V\in x^{m} L^\infty,$ for all $m>0$, implies that $e^{ m d_{g_0}(w,0)} V(w)\in L^\infty$ for all $m$, and so $V$ is super-exponentially decaying with respect to the distance function to the origin.  Therefore
\[
G_V'(r)=- \int_{-\infty}^{-\ha} 4 r  V'\bigl( T(r,t)\bigr)  \frac{1}{A(1-A^2)} \ dt.
\]
So, if $\pm V'(\rho)\geq 0$ for all $\rho$,  then $\mp G_V'(r)\geq 0$ for all $r$.  Also notice that $G_V'(r)=0$ for all $r\geq r_0$ if and only if $ V'\bigl( T(r,t)\bigr) =0$ for all $r\geq r_0$ and all $t\in (-\infty,-\ha)$.  But since $T(r,t)= d_{g_0}(w,0)$ with $w=  \theta + \frac{t\theta+\xi}{t^2+|\xi|^2}$,  this implies that $V'(\rho)=0$ for all points $w \in \mhn$ such that $d_{g_0}(w,0)\geq r_0.$ But if $V(\rho) \rightarrow 0$ as $r\rightarrow 1,$ it follows that $V(\rho)=0$ for $\rho\geq r_0.$  So we conclude that if $V'(\rho)\not=0$ if $V(r)\not=0,$ the same is true for $G_V(r).$

 So we may just assume that $V(\rho)$ nonpositive and increasing and  therefore $G_V(r)$ is nonnegative and decreasing. 
  If we have two potentials $V_1(\rho)$ and $V_2(\rho)$ for which \eqref{limit-meas} holds, then
\beq\label{MES1}
\operatorname{meas}\bigl( (G_{V_1})^{-1}(\alpha,\beta) \bigr)= \operatorname{meas}\bigl( (G_{V_2})^{-1}(\alpha,\beta) \bigr) \text{ for all } (\alpha,\beta)\subset (0,\infty).
\eeq
 In particular,  we must have $A\defn G_{V_1}(0)=G_{V_2}(0)\not=0$. Otherwise, say for example that  $G_{V_1}(0)>G_{V_2}(0)>0$,  since the functions are decreasing, there would be an interval $(\alpha,\beta)\subset (G_{V_2}(0), G_{V_1}(0)) $ that would violate \eqref{MES1}. If  $G_{V_1}(0)=G_{V_2}(0)=0$, since the functions are decreasing, $G_{V_1}(r)=G_{V_2}(r)=0$ 
 for all $r>0$.  But in that case the injectivity of the $X$-ray transform for super-exponentially decaying  potentials, see Corollary 1.3 of Chapter 3 of \cite{Hel},  implies that $V_1=V_2=0$, which certainly is not the case.

Also, if $G_{V_1}(r)=0$ for all $r\geq R,$ then 
\[
\lim_{\eps\rightarrow 0} \operatorname{meas}\bigl( G_{V_1}^{-1} (\eps, A-\eps)\bigr)= \lim_{\eps \rightarrow 0} \bigl(G_{V_1}^{-1}(A-\eps)- G_{V_1}^{-1}(\eps)\bigr)=R.
\]
Since the same would have to be true for  $G_{V_2}$, then $G_{V_2}(R)=0$ and so $G_{V_2}(r)=0$ for all $r\geq R$.  So we conclude that we either have
\[
\begin{split}
& G_{V_j}: [0, R] \longrightarrow [0,A], \;\ j=1,2, \text{ or }  \\
& G_{V_j}: [0,\infty) \longrightarrow (0,A], \;\ j=1,2.
\end{split}
\]
In either case the functions are injective and the closure of their ranges are the same.  So  we may define  
\[
\psi(r)= \bigl(G_{V_1}\bigr)^{-1}\bigl(G_{V_2} (r)).
\]
 But  \eqref{MES1} implies that if
$\alpha<\beta$ and 
\[
\begin{split}
& \alpha= G_{V_1}(a_1),  \;\ \beta= G_{V_1}(b_1),  \;\  a_1< b_1, \\
& \alpha= G_{V_2}(a_2),  \;\ \beta= G_{V_2}(b_2),  \;\  a_2< b_2,
\end{split}
\]
then   $a_2-b_2= a_1-b_1$, and so
\[
\frac{G_{V_1}(b_1)- G_{V_1}(a_1)}{ b_1-a_1}= \frac{G_{V_2}(b_2)- G_{V_2}(a_2)}{ b_2-a_2}.
\]
and we deduce  that, for every $a_j>0$, $j=1,2$, such that $G_{V_1}(a_1)=G_{V_2}(a_2)$, we have $G_{V_1}'(a_1)= G_{V_2}'(a_2)$,
and this implies that $\psi'(\alpha)=1$, but since $\alpha$ is arbitrary, $\psi'(r)=1$, but we know that  $\psi(0)=0$ and so $\psi(r)=r$

This implies that  $G_{V_1}(r)=G_{V_2}(r)$, for all $r>0$. Again,  the injectivity of the $X$-ray transform for this class of potentials  implies that $V_1=V_2$.

\end{document}